\documentclass [a4paper,twoside,12pt]{article}

\usepackage[utf8]{inputenc}
\usepackage{amsmath,amsfonts,amssymb,color}
\usepackage{vmargin,graphicx,theorem}
\usepackage[english]{babel}
\usepackage{enumerate}
\usepackage{color}
\usepackage{pst-fill,pst-grad,pst-plot,pst-eucl,pstricks-add,pst-node}
\usepackage{todonotes}
\RequirePackage[colorlinks,linkcolor=blue,citecolor=blue,urlcolor=blue]{hyperref}

\setpapersize[portrait]{A4}
\setmarginsrb{1.5cm}{1cm}{1.5cm}{2.5cm}{1.5cm}{0cm}{0.5cm}{2cm}

\newcommand{\disp}{\displaystyle}

\newtheorem{ethm}{Theorem}[section]

\newtheorem{eprop}[ethm]{Proposition}

\newtheorem{edefi}[ethm]{Definition}

\newtheorem{erem}[ethm]{Remark}

\newcommand{\proofend}{~$\rhd$}
\newcommand{\proofbegin}{~$\lhd$}

\newenvironment{eproof}
{\noindent {\emph{\textbf{Proof}}}\\\proofbegin~}
{\proofend\\}

\newcommand{\PAR}[1]{\ensuremath{{\left(#1\right)}}} 
\newcommand{\SBRA}[1]{\ensuremath{{\left[#1\right]}}} 

\renewcommand{\phi}{\varphi}

\renewcommand{\geq}{\geqslant}

\def\disp{\displaystyle}

\newcommand{\N}{\ensuremath{\mathbf{N}}}
\newcommand{\R}{\ensuremath{\mathbf{R}}}

\newcommand{\n}{\nabla}

\newcommand{\grad}{\mathrm{grad}}

\newcommand{\beq}{\begin{equation}}\newcommand{\eeq}{\end{equation}}
\begin{document}

	\title{Long-time behaviour of entropic interpolations}
	
	\author{Gauthier Clerc\thanks{Univ Lyon, Université Claude Bernard Lyon 1, CNRS UMR 5208, Institut Camille Jordan, 43 blvd. du 11 novembre 1918, F-69622 Villeurbanne cedex, France. \href{mailto:clerc@math.univ-lyon1.fr}{clerc@math.univ-lyon1.fr}},  Giovanni Conforti\thanks{D\'epartement de Math\'ematiques Appliqu\'ees, \'Ecole Polytechnique, Route de Saclay, 91128, Palaiseau Cedex, France. \href{mailto:giovanni.conforti@polytechnique.edu}{giovanni.conforti@polytechnique.edu}},  Ivan Gentil\thanks{Univ Lyon, Université Claude Bernard Lyon 1, CNRS UMR 5208, Institut Camille Jordan, 43 blvd. du 11 novembre 1918, F-69622 Villeurbanne cedex, France. \href{mailto:gentil@math.univ-lyon1.fr}{gentil@math.univ-lyon1.fr}}}
	
	\date{\today}
	\maketitle

	\abstract{In this article we investigate entropic interpolations. These measure valued curves describe the optimal solutions of the Schr\"odinger problem \cite{Sch31}, which is the problem of finding the most likely evolution of a system of independent Brownian particles conditionally to observations. It is well known that in the short time limit entropic interpolations converge to the McCann-geodesics of optimal transport. Here we focus on the long-time behaviour, proving in particular asymptotic results for the entropic cost and establishing the convergence of entropic interpolations towards the heat equation, which is the gradient flow of the entropy according to the Otto calculus interpretation. Explicit rates are also given assuming the Bakry-\'Emery curvature-dimension condition. In this respect, one of the main novelties of our work is that we are able to control the long time behavior of entropic interpolations assuming the $CD(0,n)$ condition only.}

	
	\section{Introduction}

	In two seminal papers \cite{Sch31,Sch32}, E. Schr\"odinger asked the question of finding the most likely evolution of a cloud of Brownian particles conditionally on the observation of its empirical distribution at two different times $t=0$ and $T$. In modern language, Schr\"odinger's question is translated into an entropy minimization problem under marginal constraints, known as the Schr\"odinger problem. The discovery \cite{mikami04} that the Monge-Kantorovich problem is recovered as a short time (or small noise) limit of the Schr\"odinger problem has triggered an intense research activity in the last decade. 
	Among the reasons for this renewed interest are the fact that adding an entropic penalty in the Monge-Kantorovich problem leads to major computational advantages (see for instance  \cite{COTFNT}) and the fact that the behavior of optimal solutions, called Schr\"odinger bridges or entropic interpolations, can be precisely quantified under a curvature condition. In particular, a convexity principle akin to the celebrated displacement convexity of the entropy \cite{renesse-sturm2005} holds for Schr\"odinger bridges implying a novel class of functional inequalities and the exponential convergence of entropic interpolations towards the equilibrium configuration as the time interval between observations grows larger. This last result generalizes the exponential dissipation of the entropy along the heat flow \cite{bakry-emery1985} and, in view of the stochastic control formulation of the Schr\"odinger problem, may be regarded as a turnpike theorem for Schr\"odinger bridges. Indeed, the fact that optimal curves of dynamical control problems spend of their time around equilibrium states, called turnpikes, is known in the optimal control literature as the turnpike property \cite{mckenzie1963turnpike,trelat2015turnpike}. Motivated by these results, this article aims at improving the understanding of long time behavior of entropic interpolations under a curvature-dimension condition. In particular, we aim at quantifying the role played by the dimension, so we devote a large share of our efforts to the setting $CD(0,n)$, covering in particular the case of Brownian particles in $\R^n$, which corresponds to the original Schr\"odinger problem \cite{Sch31}. Since most of the existing literature focuses on the short time regime when the Schr\"odinger problem converges towards optimal transport, much less is known for large times: in particular no asymptotic result for large times appear to be known under the $CD(0,n)$ condition. Leaving all precise statements and definitions to the main body of the article, let us give an overview of our contributions:

	\begin{itemize}
		\item We prove at Theorem~\ref{thm-cost-consqty-asym} sharp asymptotic bounds for the entropic cost $\mathcal C_{T}(\mu,\nu)$ under $CD(0,n)$ as well as for the associated energy $\mathcal E_{T}(\mu,\nu)$. 
		The entropic cost $\mathcal C_{T}(\mu,\nu)$ (see Defintion \ref{def: cost}) is the optimal value of the Schr\"odinger problem: in the large deviations interpretation of the latter, it quantifies the asymptotic probability for the cloud of independent particles to make the transition from configuration $\mu$ at time $0$ to configuration $\nu$ at time $T$. The quantity $\mathcal{E}_T(\mu,\nu)$ is called energy owing to the Otto calculus interpretation of the Schr\"odinger problem where it plays the role of the total conserved energy of a physical system. Moreover, it expresses the derivative of the cost $\mathcal{C}_T(\mu,\nu)$ w.r.t. to the time variable $T$ (see Proposition \ref{prop: envelope}). 
		In stark contrast with the results obtained under $CD(\rho,\infty)$, the cost may diverge when $T$ goes to infinity, but not faster than $\log T$ and the exponential decay of the energy does not hold, but only an algebraic one sided estimate of the order $1/T$ can be established:
		$$
		-\mathcal{E}_T(\mu,\nu) \leq \frac{2n}{T}, \quad \mathcal{C}_T(\mu,\nu) \leq \mathcal{C}_1(\mu,\nu)+2n\log(T).
		$$
		
		The sharpness of these estimates can be seen by considering Brownian particles on $\R^n$, i.e. the classical Schr\"odinger problem. Moreover, we also obtain the two-sided asymptotic estimate $|\mathcal{E}_T(\mu,\nu)| \leq C \log(T)/T$.
		
		\item We show at Theorem~\ref{thm:interpolationtoflow}, that on a fixed time window $[0,t]$, the entropic interpolation (Schr\"odinger bridge) constructed over a growing time window $[0,T]$ converges to the gradient flow (the law of a diffusion process) when $T\rightarrow+\infty$. We also establish a rate of convergence of $\log T/T$, i.e. we prove that
		
		\begin{equation*}
		W_2(\mu_t^T,P_t^*(\mu))\leq C\sqrt{\frac{\log T}{T}},
		\end{equation*}
		where $\mu^T_t$ is the entropic interpolation, $P_t^*(\mu)$ the gradient flow and $W_2$ the Wasserstein distance between probability measures.
		The $\sqrt{\log T/T}$ rate may be suboptimal as in some concrete examples we find a rate of convergence of $1/T$. This result admits a natural interpretation in terms of Schr\"odinger's thought experiment. Indeed, by ergodicity of the underlying particle system, its configurations at times $t$ and $T$ are approximately independent. Therefore what an external observer sees at time $T$ has a small influence on the particle distribution at time $t$ and particles are expected to behave almost as if no observation was made, i.e. following the gradient flow of the entropy. 
		
		\item We show at Theorem~\ref{thm:tpikeCD(0,N)} a dissipation estimate for the Fisher information $\mathcal{I}_W$ along the entropic interpolation $(\mu_{t}^T)_{t\in[0,1]}$. This estimate tells that under $CD(0,n)$ the Fisher information, calculated at time $t$ which is of the order of $T$, decays at least as fast as $1/T$\footnote{We refer to Theorem~\ref{thm:tpikeCD(0,N)} and the main body of the article for a rigorous definition of all the objects appearing in the equation below.}:
		
		$$
		\forall T>0,\theta\in(0,1), \quad \mathcal{I}_W(\mu_{\theta T}^T) \leq \frac{n}{{{2}}T \theta (1- \theta)}.
		$$
		It is worth noticing that the decay of the Fisher information at rate $1/T$ along the gradient flow is a well known fundamental consequence of the $CD(0,n)$ condition. The sharpness of the dissipation rate we establish follows from the fact that it implies a similar estimate along the gradient flow, which is known to be sharp. Besides being interesting in its own right, one may view this result as a replacement for a turnpike theorem in a context where a classical turnpike result cannot be proven. Indeed, assuming only $CD(0,n)$ is not strong enough to ensure that the associated relative entropy functional admits a minimizer among probability measures. This translates into the fact that there is no turnpike for the stochastic control formulation of the Schr\"odinger problem. However, our estimate guarantees that optimal trajectories stay in regions where the Fisher information is small.
		
	\end{itemize}
	
	Another contribution of this work is to provide alternative proofs of exponential turnpike estimates and exponential decay of the conserved quantity under the $CD(\rho,\infty)$ condition. These results have already been obtained in \cite{conforti-tamanini2019} and \cite{backhoff2019mean}. The proofs we make in this article are done in close analogy with a toy model for entropic interpolations put forward \cite{gentil-leonard2020} and are therefore simpler to read and amenable to generalizations beyond the framework considered in the above mentioned references.
	\paragraph{Organization} In Section~\ref{sec:setting} we introduce curvature-dimension conditions, recall some basic facts about the Schr\"odinger problem and state our main hypothesis. In section~\ref{sec:toy} we prove the main results of the paper for a toy model introduced in~\cite{gentil-leonard2020}. In Section~\ref{sec:main} we lift our results from the simple setting of the toy model to the general Schr\"odinger problem. Along the way, we illustrate the sharpness (or not) of our results by means of examples. The case of the Euclidean heat semigroup is studied in more detail at Section~\ref{sec:heat}.

	\section{Setting of our work}\label{sec:setting}
	
	\subsubsection*{Markov semigroups and the curvature-dimension condition $CD(\rho,n)$}
	
	Let $(N,\mathfrak g)$ be a smooth, complete and connected Riemannian manifold. We consider the generator $L=\Delta-\nabla W\cdot \nabla$ where $\Delta$ is the Laplace-Beltrami operator, $\nabla$ is the gradient operator. $\nabla\cdot$ is the divergence operator (in order to have $\Delta=\nabla\!\cdot\nabla$) and $W:N\longrightarrow \R$ is a smooth function. The carré du champ operator $\Gamma$ is defined for any smooth functions $f,g$ by 
	$$
	\Gamma (f,g)=\frac12\PAR{L(fg)-fLg-gLf}.
	$$
	Under the current hypothesis  $\Gamma(f)=|\nabla f|^2$, which is the length of $\nabla f$ with respect to the metric $\mathfrak g$ (for simplicity, we omit  the dependence with respect to the metric). As usual, we adopt the shorthand notation $\Gamma (f)$ for $\Gamma(f,f)$. 
	The measure $d Vol$ denotes the Riemannian volume measure. 
	Whenever $Z=\int e^{-W}d Vol<+\infty$ then we set  $dm=\frac{e^{-W}}{Z}d Vol$, the corresponding probability distribution, that is reversible for $L$. When $Z=\infty$, then we set $dm=e^{-W}d Vol$: $m$ has infinite mass and is still reversible for $L$. We denote by $\mathcal P(N)$ (resp. $\mathcal P_2(N)$ and $\mathcal M(N)$), the set of probability measures  on $N$ (resp. probability measures admitting a second moment and the set of positive measures).   We assume that $L$ is the infinitesimal generator of a Markov semigroup in the sense proposed in~\cite[Sec.~3.2]{bgl-book}, that is to say, $(N, \Gamma, m)$ is a full Markov triple. The Markov semigroup is denoted $(P_t)_{t\geq0}$, and is identified with the map $(t,x)\mapsto P_tf(x)$ solution of the parabolic equation  
	\begin{equation}
	\label{eq-45}
	\left\{
	\begin{array}{l}
	\partial_tu=Lu\\
	u(0,\cdot)=f(\cdot), 
	\end{array}
	\right.
	\end{equation}
	for function $f\in {\rm L}^2(m)$.  The Markov semigroup admits a Markov kernel $p_t(x,dy)$  with density $p_t(x,y)$ against the invariant measure $m$, that is for all functions $f\in {\rm L}^2(m)$
	$$
	P_tf(x)=\int f(y)p_t(x,dy)= \int f(y) p_t(x,y)dm(y).
	$$
	We also introduce the dual semigroup $(P_t^*)_{t\geq 0}$ acting on absolutely continuous probability measures $\mu\in\mathcal P(N)$ as follows
	\begin{equation}
	\label{eq-76}
	P_t^*(\mu)=P_t\PAR{e^{W}\frac{d\mu}{d Vol}}m=P_t \PAR{\frac{d\mu}{dm}}m \in \mathcal{P}(N).
	\end{equation}
	{One finds that $(t,x)\mapsto  \frac{dP_t^*(\mu)}{dVol}$ is a solution  of the  Fokker-Planck equation, 
		\begin{equation}
		\label{eq-101}
		\partial_t u_t=L^*u_t=\Delta u_t+\nabla\cdot(u_t \nabla W)=\nabla\cdot(u_t\nabla(\log u_t+W)),
		\end{equation}
		starting from $\frac{d \mu}{dVol}$. 
	}
	
	Following the seminal work of Bakry-\'Emery~\cite{bakry-emery1985}, we say that the semigroup satisfies the curvature-dimension condition $CD(\rho,n)$  with $\rho\in\R$ and $n\in(0,\infty]$ if for any smooth function $f$ defined on $N$, 
	\begin{equation}
	\label{eq-52}
	\Gamma_2(f)\geq \rho\Gamma(f)+\frac{1}{n}(Lf)^2,
	\end{equation}
	where $\Gamma_2(f)=\frac{1}{2}L\Gamma(f)-\Gamma(f,Lf)$ is the iterated carré du champ operator. Following again~\cite{bakry-emery1985}, the curvature-dimension condition $CD(\rho,n)$ with $\rho\in\R$ and $n\geq d$ (where $d\in \N^*$ is the dimension of $N$) is equivalent to the following inequality on tensors
	$$
	{\rm Ric }(L):={\rm Ric}_{\mathfrak g}+\nabla\nabla W-\rho\mathfrak g\geq\frac{1}{n-d}\nabla W\otimes\nabla W.
	$$
	When $n=d$, then we need to impose $W=0$ in the above. In particular, if ${\rm Ric}_{\mathfrak g}\geq \rho\mathfrak g$ for some $\rho\in\R$, then the Laplace-Beltrami operator $\Delta$ satisfies the $CD(\rho,d)$ condition. On $\R^n$, if $W(x)=|x|^2/2$, then $L$ satisfies $CD(1,\infty)$, whereas if $W=0$, $L$ is the Euclidean Laplace operator  which satisfies $CD(0,n)$.

	\subsubsection*{Statement of the Schrödinger problem}
	
	In this section, we recall some basic facts about the Schrödinger problem following the presentation of~\cite{leonard2014}, see also~\cite{follmer1988}. In order to do so, we need to introduce the relative entropy functional, defined for any probability measure $\mathbf q$ and a positive measure $\mathbf r$ on the same measurable space as follows
	\begin{equation}
	\label{eq-41}
	H(\mathbf q|\mathbf r)=
	\left\{
	\begin{array}{ll}
	\disp\int \log \frac{{d}\mathbf q}{{d}\mathbf r} d\mathbf q\in(-\infty,+\infty],& {\rm if} \,\,\mathbf q\ll \mathbf r;\\
	\disp+\infty &otherwise.
	\end{array}
	\right.
	\end{equation}
	{This definition is meaningful when $\mathbf r$ is a probability measure but not necessarily when $\mathbf r$ is unbounded. Assuming that $r$ is $\sigma$-finite, then there exists a function $W: M \rightarrow [1, \infty)$ such that $z_W:= \int e^{-W}d \mathbf r < \infty$. Hence we can define a probability measure $\mathbf r_W:=z_W^{-1}e^{-W}\mathbf r$ and for every measure $\mathbf q$ such that $\int W d\mathbf q < \infty$
		$$
		H(\mathbf q|\mathbf r):=  H(\mathbf q |r_W)- \int W d\mathbf q - \log(z_W),
		$$
		where $H(\mathbf q|\mathbf r_W)$ is defined by the equation~(\ref{eq-41}). Hence to ensure the existence and finiteness of $H(\mathbf q| \mathbf r)$ it is enough to assume that $H(\mathbf q|\mathbf r_W) < \infty$ and $W \in L^1(\mathbf q)$.} For a given $T>0$,  let $\Omega=\mathcal C([0,T],N)$ be the set of continuous paths from $[0,T]$ to $N$ on which we define the probability measure $R_x\in\mathcal P(\Omega)$ as the law of a Markov process with generator $L$, started at  $x$.  Finally we define the positive measure $R^T(\cdot)$ by
	$$
	R^T(\cdot)=\int R_x(\cdot)dm(x)\in\mathcal M(\Omega).
	$$
	
	For a given pair of probability measures $\mu,\nu\in\mathcal P(N)$, the Schrödinger problem is
	\begin{equation}
	\label{eq-40}
	{\rm Sch}_T(\mu,\nu)=\inf \{H(Q|R^T),\,\,Q\in\mathcal P(\Omega),\,\,Q_0=\mu,\,\,Q_T=\nu\},
	\end{equation}
	where $Q_0=X_0\# Q$ and $Q_T=X_T\# Q$. Here $(X_t)_{t\in[0,T]}$ is the canonical process and $X_0\# Q\in\mathcal P(N)$ is the image measure of  $Q$ by  $X_0$, that is, for any test function $h$,  $\int hd X_0\# Q=\int h(X_0)dQ(X)$. In other words, it is the problem of minimizing the relative entropy $H(\cdot|R^T)$ among all path probability measures $Q\in\mathcal P(\Omega)$ with prescribed initial marginal $\mu$ and final marginal $\nu$, that is $X_0\# Q=\mu$ and $X_T\# Q=\nu$. 
	
	Also, notice that the Schrödinger problem admits a static formulation, that is 
	$$
	Sch_T(\mu,\nu)=\inf \left\{ H(\gamma |  R_{0T}^T), \ \gamma \in \mathcal{P}(N \times N), \gamma_0=\mu, \, \gamma_1=\nu\right\},
	$$
	where $R_{0T}^T=(X_0,X_T) \# R^T$ is the joint measure of initial and final position of $R^T$, see~\cite{leonard2014}.
	
	\subsubsection*{Fundamental results on the Schrödinger problem and usual hypothesis}

	In order to ensure the existence of an optimal solution we suppose throughout this article that  for any $T>0$, there exist two non negative measurable functions $A,B$ such that
	\begin{enumerate}[(i)]
		\item $p_T(x,y)  \geq e^{-A(x)-A(y)}$,  $\forall x,y\in N$;
		\item $\int e^{-B(x)-B(y)}p_T(x,y)\,m(dx) m(dy)< \infty$;
		\item $\int (A+B)\, d\mu, \int (A+B)\, d\nu<\infty$;
		\item The quantities $H(\mu|m)$ and $H(\nu|m)$ are well defined as explained above, and are finite.
	\end{enumerate}
	
	Let us notice that hypothesis (i) and (ii) are satisfied for a large class of Markov semigroup, in particular for the one studied in this paper, semigroup satisfies a $CD(\rho,\infty)$ conditions with $\rho\in\R$. For more details, we refer to~\cite{leonard2014,gigli2018second}.
	Under these assumptions, it is proven at~\cite[Theorem~2.12]{leonard2014} that the entropic cost  ${\rm Sch}_T(\mu,\nu)$ is finite and has a unique minimizer $Q\in\mathcal P(\Omega)$. Moreover the minimizer has the following product form
	\begin{equation}
	\label{eq-42}
	\frac{dQ}{dR^T}=f(X_0)g(X_T),
	\end{equation}
	for some measurable and positive functions $f$ and $g$ on $N$. The above formula implies that if we denote by $(\mu_t^T)_{t\in[0,T]}$ the \emph{entropic interpolation}
	$$
	\mu_t^T=X_t\#Q=Q_t\in\mathcal P (N), \,\,t\in[0,T], 
	$$ 
	then we have 
	\begin{equation}
	\label{eq-43}
	\mu_t^T=P_tfP_{T-t}g\,m.
	\end{equation}
	All of these results can be found for instance in the survey~\cite{leonard2014}.

	\subsubsection*{The Benamou-Brenier-Schrödinger minimization problem}

	In analogy with the Benamou-Brenier fluid dynamic formulation of optimal transport \cite{benamou-brenier2000}, we can recast the Schrödinger problem as a minimization problem among absolutely continuous curves on $\mathcal P_2(N)$ with respect to the Wasserstein distance.
	We recall that the Wasserstein distance is defined as follow, for every $\mu,\nu\in\mathcal P_2(N)$, 
	$$
	W_2(\mu,\nu)=\inf \sqrt{\iint d(x,y)^2d\pi(x,y)},
	$$
	where the infimum is running over all $\pi\in\mathcal P (N\times N)$ with marginals $\mu$ and $\nu$.

	Following the presentation of~\cite[Chap.~1]{ambrosio-gigli2008}, we recall that a path $[0,T]\ni t\mapsto \mu_t\in\mathcal P_2(N)$ is
	absolutely continuous if there exists a non negative function $l\in L^2([0,T])$ such that for any $0\leq s\leq t\leq T$, 
	$$
	W_2(\mu_t,\mu_s)\leq\int_s^t l(r)dr.
	$$
	In that case, one can define the metric derivative $(|\mu'|(t))_{t\in[0,T]}$, as 
	$$
	|\mu'|(t):=\underset{s\rightarrow t}{\limsup}\frac{W_2(\mu_t,\mu_s)}{|t-s|}\in L^1([0,T]),
	$$ 
	a.e.  in $[0,T]$,~\cite[Thm~1.1.2]{ambrosio-gigli2008}.
	
	Following~\cite{erbar}, for any absolutely continuous path $(\mu_t)_{t\in[0,T]}$, there exists a unique vector field $(t,x) \mapsto V_t(x)$  such that, a.e. in $[0,T]$, $\int|V_t|^2d\mu_t<+\infty$ (where $|V_t|$ is the length of $V_t$ with respect to metric $\mathfrak g$) and satisfying in a weak sens, the continuity equation
	\begin{equation}\label{eq-112}
	\partial_t\mu_t+\nabla\cdot(\mu_t V_t)=0,
	\end{equation}
	and, a.e. in $[0,T]$, 
	$$
	|V_t|_{L^2(\mu_t)}=|\mu'|(t).
	$$

	The vector field $V_t$ is in fact  a limit in $L^2(\mu_t)$ of gradient of smooth compactly supported functions in $N$, as it is explained in the section related to the Otto calculus, cf. page~\pageref{sec-toto}. For every $t \in [0,T]$ we denote 
	\begin{equation}\label{eq-130}
	\dot \mu_t:=V_t,
	\end{equation}
	and we call $\dot \mu_t$ the velocity of the path $(\mu_t)_{t \in [0,1]}$ at time $t$.

	For instance, in the case of the generalized Fokker-Planck equation~\eqref{eq-101}, the velocity of the path $(\nu_t)_{t\geq0}$ is 
	\begin{equation}
	\label{eq-102}
	\dot\nu_t=-\nabla\left(\log \frac{d\nu_t}{dVol}+W\right)=-\nabla\left(\log \frac{d\nu_t}{dm}\right). 
	\end{equation}

	For every $\mu \in \mathcal P_2(N)$ and any vector field $V$ and $W$  in $L^2(\mu)$, we denote by $\langle V , W \rangle_{\mu}=\int V\!\cdot\! W d\mu $ the natural scalar product of $L^2(\mu)$ and $|V|_{\mu}$ the associated norm.
	
	\begin{edefi}[Entropic cost function]\label{def: cost}
		For any measures $\mu,\nu\in\mathcal P_2(N)$, let define 
		\begin{equation}
		\label{eq-47}
		\mathcal{C}_T(\mu,\nu)=\inf\left\{\int_0^T\SBRA{{|\dot\mu_s|_{\mu_s}^2}+{\mathcal I_{ W}(\mu_s)}}ds\right\}\in[0,\infty],
		\end{equation}
		where the infimum runs over all absolutely continuous paths $(\mu_s)_{s\in[0,T]}$ satisfying $\mu_0=\mu$ and $\mu_T=\nu$.  In the above, for any probability measure $\mu\in\mathcal P_2(N)$, $\mathcal I_W$ denotes the Fisher information,
		\begin{equation}
		\label{eq-69}
		\mathcal I_W(\mu)=
		\displaystyle\int\Gamma\PAR{\log \frac{d\mu}{dm}}d\mu= \int \Gamma \PAR{\log\left( \frac{d\mu}{d Vol}\right)+W}d\mu\in[0,+\infty],
		\end{equation}
		if quantities are well defined (smooth enough for instance) and $+\infty$ otherwise.
	\end{edefi}
	Through a simple change of variable, if we define
	$$
	\mathcal A_T(\mu,\nu)=\inf\left\{\int_0^1\SBRA{{|\dot\mu_s|_{\mu_s}^2}+T^2{\mathcal I_W(\mu_s)}}ds\right\},
	$$
	where now the infimum is running over all paths $(\mu_s)_{s\in[0,1]}$ absolutely continuous with respect to the Wasserstein distance, satisfying the condition $\mu_0=\mu$ and $\mu_1=\nu$,   then we have 
	$$
	\mathcal A_T(\mu,\nu)=T \mathcal{C}_T(\mu,\nu).
	$$
	Let us notice that, for simplicity, the definition of the cost $\mathcal A_T$ differs by a factor 2 from the one defined in~\cite{gentil-leonard2020}. For any probability measure $\mu\in\mathcal P(N)$ we define the relative entropy functional
	\begin{equation}
	\label{eq-100}
	\mathcal{F}(\mu)=H(\mu |m)
	\end{equation}
	The following result relates precisely all the variational problems encountered so far. 
	
	\begin{ethm}[Benamou-Brenier-Schrödinger formulation]
		\label{thm-1}
		For any  compactly supported measures $\mu,\nu\in\mathcal P(N)$
		\begin{equation}
		\label{eq-46}
		\mathrm{Sch}_T(\mu,\nu)=\frac{\mathcal A_T(\mu,\nu)}{4T}+\frac{1}{2}\PAR{\mathcal{ F}(\mu)+\mathcal{ F}(\nu)}=\frac{ \mathcal{C}_T(\mu,\nu)}{4}+\frac{1}{2}\PAR{\mathcal{ F}(\mu)+\mathcal{ F}(\nu)}.
		\end{equation}
	\end{ethm}
	Versions of this result have been proven in different papers~\cite{chengeorgiou2016, gentil-leonard2017, gentil-leonard2020, gigli-tamanini2020,clerc}. 
	
	\subsubsection*{Otto calculus, Hessian of $\mathcal F$ and Newton equation}
	\label{sec-toto}

	Otto calculus, developed in the seminal papers~\cite{jko1998,otto2001, otto-villani2000}, allows to formally view the space $\mathcal{P}(N)$ as an infinite dimensional Riemannian manifold. This viewpoint has already proven to be extremely useful as it provides an interpretation of a large class of dissipative PDEs as gradient flows, greatly facilitating the task of obtaining entropy dissipation estimates if the entropy under consideration is displacement convex. In this short section, we give a very concise introduction to Otto calculus, explaining at the formal level why, although entropic interpolations are not gradient flows, adopting such viewpoint still gives precious insights. Our presentation is based on~\cite{gentil-leonard2020}, to which we refer for more details. In this article, we use Otto calculus as an heuristic guideline. However, many of the following statement can be turned into rigorous statements, see the monograph \cite{gigli2012,erbar}.

	Heuristically, the tangent space at $\mu\in\mathcal P_2(N)$ is identified with
	\[T_{\mu} \mathcal P_2(N)=\overline{\{\nabla \varphi,\,\,\phi \in C_c^{\infty}(N)\}}^{L^2(\mu)}.\]
	The Riemannian metric on $T_{\mu} \mathcal P_2(N)$ is then defined via the scalar product $L^2(\mu)$ that we introduced before and denoted $\langle\cdot,\cdot\rangle_\mu$. Such metric is often referred to the Otto metric and it can be seen that the geodesics associated to the Otto metric are the displacement interpolations of optimal transport. Using this, a straightforward computation implies that the gradient of the entropy $\mathcal{F}$ at $\mu$ is given by
	$$
	\grad_{\mu} \mathcal{F}= \nabla \log \left( \frac{d \mu}{dm}\right) \in T_{\mu} \mathcal P_2(N).
	$$
	Accordingly, we can rewrite the Fisher information functional $\mathcal{I}_{W}$ as
	$$
	\mathcal{I}_W(\mu):=|\grad_{\mu} \mathcal F|^2_{\mu}=:\mathbf{\Gamma}(\mathcal{F})(\mu),
	$$
	where $\mathbf{\Gamma}(\mathcal{F})$ can be interpreted as the carré du champ operator applied to the functional $\mathcal{F}$.
	In light of~\eqref{eq-102}, we can now view the semigroup $(P_t^*)_{t\geq 0}$ as the gradient flow of the function $\mathcal{F}$, that is to say 
	$$
	\dot\nu_t=-\grad_{\nu_t} \mathcal{F}.
	$$
	Now, we turn our attention to the second order Otto calculus introducing covariant derivatives and Hessians. A remarkable fact is that the Hessian of the entropy functional $\mathcal{F}$ can be expressed in terms of the $\Gamma_2$ operator. Indeed we have (see for instance~\cite[Sec 3.3]{gentil-leonard2020}) 
	$$
	\forall\,\mu \in \mathcal P_2(N),\, \nabla \phi, \nabla \psi \in T_{\mu} \mathcal{P}_2(N), \quad\mathrm{Hess}_{\mu} \mathcal{F}(\nabla\phi , \nabla\psi)= \int \Gamma_2(\phi, \psi) d \mu.
	$$
	
	At this point we can see that the curvature-dimension condition $CD(\rho,n)$ ($\rho\in\R$, $n>0$) is equivalent to the differential inequality
	\begin{equation}
	\label{eq-103}
	\forall \mu \in \mathcal{P}_2(N), \, \nabla\phi\in T_{\mu} \mathcal P_2(N), \quad \mathrm{Hess}_{\mu} \mathcal{F}(\nabla\phi, \nabla\phi) \geq \rho |\nabla\phi|_{\mu}^2+ \frac{1}{n} \langle \grad_{\mu} \mathcal F, \nabla\phi \rangle_{\mu}^2.
	\end{equation}
	From the work of~\cite{erbar-kuwada2015}, we know that  the infinitesimal generator $L$ satisfies the curvature-dimension condition~\eqref{eq-52} if and only if the functional $\mathcal F$ satisfies the differential equation~\eqref{eq-103}. A crucial fact about entropic interpolations, i.e. the optimizers of~\eqref{eq-47}, is that they solve a second order differential equation. In order to state the equation, we need to introduce the notion of acceleration of a flow $(\mu_t)_{t\in[0,T]}$. As in a  finite dimensional Riemannian manifold, the acceleration of a curve is defined as the covariant derivative of the velocity field along the curve itself. Recalling the definition of velocity $\dot{\mu}_t$ we gave through \eqref{eq-112}, it turns out that the acceleration, which we denote $\ddot{\mu}_t$ is given by 
	\begin{equation}
	\label{1580}
	\ddot{\mu}_t= \nabla \left( \frac{d}{dt} \varphi_t+ \frac{1}{2} |\nabla \varphi_t|^2\right) \in T_{\mu_t}\mathcal{P}_2(N),
	\end{equation}
	in the case where the velocity is given by $\dot \mu_t^T=V_t=\nabla\varphi_t$ (if the velocity is not the gradient of a function, the expression of the acceleration is less pleasant). It has been noted in~\cite[Theorem 1.2]{conforti2017} (see also~\cite[Sec 3.3. and Proposition~3.5]{gentil-leonard2020}) that the entropic interpolation $\left( \mu_t^T \right)_{t \in [0,T]}$  is a solution of the following second order equation
	
	\begin{equation} \label{newteq}
	\ddot \mu_t^T= \frac{1}{2} \grad_{\mu_t^T}\mathbf{\Gamma}(\mathcal{F})=\mathrm{Hess}_{\mu_t^T}\mathcal{F}\big(\grad_{\mu_t^T}\mathcal{F}\big) \, \in T_{\mu^T_t} \mathcal P_2(N).
	\end{equation}

	Let us mention that formulas~\eqref{eq-103},~\eqref{1580} and~\eqref{newteq} can be justified rigorously but actually we only use it  as an heuristic guideline. 
	
	We call the above a Newton equation, in analogy with Newton's law $\ddot{X}= F(X)$, which describes the evolution of a particle in a force field. In the rest of the paper, we shall heavily exploits this analogy in order to obtain the main results.
	
	\medskip

	\section{The finite dimensional case}\label{sec:toy}
	\label{sec-3}
	In this section we study a toy model introduced in~\cite[Sec.~2]{gentil-leonard2020}. Despite its simplicity, this model already captures quite well the geometric structure of the Schr\"odinger problem. In fact, we shall see in the next section that the results obtained for the toy model transfer with little effort to the Schr\"odinger problem. Let $F: \R^n \mapsto \R$ be a twice differentiable function with $d>0$. We note $F'$ (resp. $F''$) the gradient (resp. the Hessian) of $F$. For every $T>0$ and $x,y \in \R^n$, the toy model is the following optimization problem 
	\begin{equation}
	\label{eq-48}
	C_T(x,y)=\inf \left\{ \int_0^T \SBRA{{|\dot \omega_s|^2}+{|F'(\omega_s)|^2}}dt\right\},
	\end{equation}
	where the infimum taken over all smooth  paths from $[0,T]$ to $\R^n$ such that $\omega_0=x$ and $\omega_T=y$ and $\dot\omega_s=\frac{d}{ds}\omega_s$. A standard variational argument shows that any minimizer  $\left( X_t^T \right)_{t \in [0,T]}$  of~\eqref{eq-48} satisfies Newton's system
	\begin{equation}
	\label{geq-1}
	\left\{
	\begin{array}{l}
	\disp\ddot X_t^T =\frac{1}{2} (|F'|^2)'(X_t^T)=F''(X_t^T)F'(X_t^T),\\
	\disp X_0^T=x,\,\,X_T^T=y,
	\end{array}
	\right.
	\end{equation}
	and is called an $F$-interpolation between $x$ and $y$. If $(X_t^T)_{t\in[0,1]}$ is an $F$-interpolation, then from Newton's equation~\eqref{geq-1} we get that the quantity 
	$$
	{E}_T(x,y)=\big|\dot X_t^T\big|^2-\big|F'(X_t^T)\big|^2,
	$$
	is conserved, i.e. it does not depend on $t$. 
	Let $(S_t)_{t \geq 0}$ be the gradient flow semigroup of $F$ that is for every $x \in \R^n$, $(S_t(x))_{t \geq 0}$ is the only solution of
	\begin{equation} \label{eq-1001}
	\left\{ \begin{array}{l}
	\frac{d}{dt} S_t(x)=-F'(S_t(x)), \ t \geq 0 \\
	S_0(x)=x.
	\end{array}\right.
	\end{equation}
	Heuristically, the best way to minimize $C_{T}(x,y)$ is to follow closely the gradient flow for most of the time, and only when final time $T$ is very close, depart from it to reach the target destination $y$. In terms of Schr\"odinger's thought experiment, this means that the effect of the observation made at $T$ affects only slightly the dynamics of the particle systems at time $t$, provided $T-t$ is large. Using the language of control theory, what we are saying is that $F$-interpolation satisfy the turnpike property \cite{trelat2015turnpike}.
	This leads to believe that, for $t\geq0$, 
	$$
	X_t^T \underset{T \rightarrow \infty}{\rightarrow} S_t(x). 
	$$
	In Sections~\ref{sec-3.2} and~\ref{sec-3.3} we establish a quantitative form of this convergence results under two different types of convexity hypothesis on $F$, which are finite dimensional analogs of the $CD(\rho,+\infty)$ and $CD(0,n)$ conditions. Indeed, inspired by \eqref{eq-103},
	we say that $F$ is $(\rho,n)$-convex for some $\rho\in\R$ and $n\in(0,\infty]$ if
	$$
	F''\geq \rho{\rm Id}+\frac{1}{n}F' \otimes F'.
	$$
	Here, we only treat the case where $F$ is $(\rho,\infty)$ or $(0,n)$-convex.

	\subsection{Two examples in finite dimension}
	\label{sec-3.1}
	
	To build intuition, we start working on two examples, which allow for explicit calculations. In both cases, we provide precise estimates for the three quantities of interest (calculations are detailed in Appendix~\ref{AppendiceA}):
	\begin{itemize}
		\item  the cost $C_T(x,y)$; 
		\item the conserved quantity ${E}_T(x,y)$;
		\item the distance between the $F$-interpolation and the gradient flow $\left|X_t^T-S_t(x)  \right|$.
	\end{itemize}

	\subsubsection{A $(1,\infty)$-convex function}
	\label{sec-3.2.0}
	
	Let consider  $F(x)=|x|^2/2$, $x\in\R^n$.  Then $F''={\rm Id}$, that $F$ is $(1,\infty)$-convex. We find that  
	\begin{itemize}
		\item The gradient flow starting from $x\in\R^n$, is given by
		$
		S_t(x)=e^{-t}x,\,\,t\geq0.
		$
		
		\item  The $F$-interpolation $\left(X_t^T \right)_{t \in [0,T]}$ between $x$ and $y$ is given by
		$
		X_t^T=S_t\left(\alpha_T\right)+S_{T-t}\left(\beta_T\right)
		$, for $t\in[0,T]$ where  
		$\alpha_T=\frac{x-ye^{-T}}{1-e^{-2T}}$ and $\beta_T=\frac{y-xe^{-T}}{1-e^{-2T}}$.
		
		\item  For all $x,y \in \R^n$ and $T>0$, the conserved quantity is given by
		$
		{E}_T(x,y)\!=\!-4e^{-T}\alpha_T  \beta_T 
		$
		and there exists a constant $C>0$ (depending on $x,y$) such that
		$$
		\left|E_T(x,y) \right| \leq e^{-T}C,\,\,\, T>0.
		$$

		\item  The cost function is given by
		$
		C_T(x,y)=\left(1-e^{-2T} \right)\left(|\alpha_T|^2 +|\beta_T|^2\right),
		$
		$x,y \in \R^n$, 
		and therefore there exists a constant $C>0$ (depending only on $x$ and $y$) such that for any $T>0$,
		$$
		C_T(x,y) \leq C.
		$$
		\item For all $x,y \in \R^n$, the distance between entropic interpolation and gradient flow is given by
		$$
		\left|X_t^T - S_t(x) \right|=\frac{2\sinh(t)}{1-e^{-2T}}\left|y-xe^{-T} \right|e^{-T},\,\,T>0.
		$$

	\end{itemize}
	
	As a conclusion in this example, the $F$-interpolation converges exponentially fast toward the gradient flow. 
	\subsubsection{A $(0,1)$-convex function}
	\label{sec-3.2.1}
	
	Let $F(x)=-\log(x)$, for any $x>0$.  Since $F''=(F')^2$, then $F$ is a $(0,1)$-convex function. All computations are explained in Appendix~\ref{AppendiceA2}.
	\begin{itemize}
		\item The gradient flow from $x>0$, is given by
		$
		S_t(x)=\sqrt{2t+x^2},\,\,\, t\geq0.
		$
		\item  For all $x>0$ and $T>0$ the conserved quantity is given by
		$
		{E}_T(x,x)=\frac{-x^2-\sqrt{x^4+T^2}}{T^2/2},
		$
		and there exists a constant $C>0$ (depending on $x$) such that
		$$
		\left|{E}_T(x,x) \right| \leq \frac{C}{T},\,\,\,T>0.
		$$
		\item The cost function satisfies 
		$
		C_T(x,x) \underset{T \rightarrow \infty}{\sim} 2 \log(T).
		$
		\item  The $F$-interpolation between $x$ and $x$ is given by
		$$
		\forall t \in [0,T], \ X_t^T= \sqrt{x^2+t^2{E}_T(x,x)+2t\sqrt{1+{E}_T(x,x)x^2}}.
		$$
		\item There exists a constant $C>0$ (depending on $x$) such that
		$$
		|X_t^T-S_t(x)| \leq \frac{C}{T},\,\,\, t\in[0,1], \,\,\,T>0.
		$$
		
	\end{itemize}
	\subsection{The  $(\rho,\infty)$-convex case}
	\label{sec-3.2}

	In this section we are assuming that $F$ is a smooth and positive  $\rho$-convex function for some $\rho>0$, that is 
	\begin{equation}
	\label{eq-50}
	F'' \geq \rho \,{\rm Id}.
	\end{equation}
	Since $\rho>0$,  there exists $x^*\in\R^n$ such that $\inf F=F(x^*)$. 
	
	Under this convexity condition, the cost is bounded, that is for all $x,y \in \R^n$ and $T>0$
	\begin{equation}
	\label{eq-49}
	C_T(x,y) \leq  2 \frac{1+e^{- \rho T}}{1-e^{-\rho T}}(F(x)+F(y)-2F(x^*))
	\end{equation}
	see~\cite[Cor~2.13]{gentil-leonard2020}. The above result can be reinforced  as follows, with the same proof,
	\begin{equation}\label{eq-201}
	C_T(x,y) \leq   \inf_{t\in(0,T)}\left\{ 2 \frac{1+e^{- 2\rho t}}{1-e^{-2\rho t}}(F(x)-F(x^*))+2 \frac{1+e^{- 2\rho (T-t)}}{1-e^{-2\rho (T-t)}}(F(y)-F(x^*))\right\}.
	\end{equation}
	Thus, the cost is bounded by a constant, depending only on $x$ and $y$. To quantify how far the $F$-interpolation $(X^T_t)_{t\in[0,T]}$  is from the gradient flow we introduce the function
	$$
	\phi_t^T:=F'(X_t^T)+\dot X_t^T, \,\, t \in [0,T]. 
	$$
	First we  control the vector field  $(\phi_t^T)_{t\in[0,T]}$. 
	\begin{eprop} 
		\label{ggprop-2}
		For all $x,y\in\R^n$, $T>0$ and $t\in(0,T)$ we have 
		$$
		|\phi_t^T|^2 \leq  \frac{2\rho}{\exp(2\rho(T-t))-1}(\mathcal{C}_T(x,y)+2F(y)-2F(x)) .
		$$
		In particular,
		\begin{equation}\label{eq-206}
		|\phi_t^T|^2 \leq  \frac{8\rho}{\exp(2\rho(T-t))-1}\Big(\frac{e^{- 2\rho t}}{1-e^{-2\rho t}}(F(x)-F(x^*))+ \frac{1}{1-e^{- 2\rho (T-t)}}(F(y)-F(x^*)) \Big).
		\end{equation}
	\end{eprop}

	\begin{eproof}
		Newton equation~\eqref{geq-1} implies that $
		\frac{d}{dt}\phi_t^T=F''(X_t^T)\phi_t^T$. Combining with~\eqref{eq-50} we get
		$$
		\frac{d}{dt} |\phi_t^T|^2 \geq 2\rho |\phi_t^T|^2.
		$$
		Therefore for all $t\leq s\leq T$ we find
		$|\phi_s^T|^2 \geq \exp(2\rho(s-t))|\phi_t^T|^2,
		$ 
		and integrating this bound over $[t,T]$ we get
		\begin{equation*}
		\int_t^T|\phi_s^T|^2 ds \geq \frac{\exp(2\rho(T-t))-1}{2\rho}|\phi_t^T|^2.
		\end{equation*}
		Observing that $\int_t^T|\phi_s^T|^2 ds \leq \mathcal{C}_T(x,y)+2F(y)-2F(x)$ and using \eqref{eq-201} we obtain the desired results.
	\end{eproof}

	\begin{ethm}[Convergence of the $F$-interpolation]
		\label{gthm-1}
		For all $x,y\in\R^n$, $T>0$ and $t \in (0,T)$
		$$
		|X_t^T-S_t(x)|\leq t \exp (- \rho T)  \sqrt{\frac{2 \rho}{\exp(-2\rho t)-\exp(-2\rho T)}\left(C_T(x,y)+2F(y)-2F(x)\right)}.
		$$
		Furthermore there exists a constant $C>0$ depending only on $x$ and $y$ such that for every $t \geq 0$ and $T>t$
		$$
		|X_t^T-S_t(x)| \leq C \frac{t \exp(-\rho T)}{\sqrt{\exp(-2\rho t)-\exp(-2\rho T)}}.
		$$

	\end{ethm}
	
	\begin{eproof}
		Let $0 \leq t \leq T-1$. Whence by the Cauchy-Schwarz inequality and the Proposition~\ref{ggprop-2} we have
		\begin{equation*}
		\begin{split}
		\frac{d}{dt}\frac{|X_t^T-S_t(x)|^2}{2}&=\langle\dot X_t^T+F'(S_t(x)), X_t^T-S_t(x) \rangle \\
		&=-\langle F'(X_t^T)-F'(S_t(x)), X_t^T-S_t(x) \rangle+\langle \dot X_t^T +F'(X_t^T),X_t^T-S_t(x) \rangle \\
		&\leq |\varphi_t^T||X_t^T-S_t(x)|. \\
		\end{split}
		\end{equation*}
		The result follow from integration of this inequality and the Proposition~\ref{ggprop-2}. 
	\end{eproof}
	
	According to the example given in Section~\ref{sec-3.2.0}, Theorem~\ref{gthm-1} gives the optimal rate for the convergence.

	\subsubsection{Turnpike property}
	
	Under the hypothesis that $F$ is $\rho$-convex with $\rho>0$ as defined in~\eqref{eq-50}, it is well known that the gradient flow $S_t$ dissipates $F$ at exponential rate $2\rho$. This mean that,
	\begin{equation*}
	F(S_{T}(x))-F(x^*) \leq \exp(-2\rho T)(F(x)-F(x^*)). 
	\end{equation*}
	The aim of this subsection is to show that a similar estimate holds replacing the gradient flow with the $F$-interpolation. A fundamental ingredient needed for the proof of this result is the following exponential upper bound for the conserved quantity $E_{T}(x,y)$.
	
	\begin{eprop}
		For all $x,y\in \R^n$, $T>0$ 
		
		\begin{equation}\label{eq-202}
		|E_{T}(x,y)|\leq \frac{2\rho}{\exp(\rho T)-1} \sqrt{{C}^2_T(x,y)-4(F(x)-F(y))^2}
		\end{equation}
		
	\end{eprop}
	
	\begin{eproof}
		Denoting by $\langle\cdot ,\cdot \rangle$ the inner product in $\R^n$, we obtain
		
		\begin{equation*}
		|E_{T}(x,y)|= |\langle \dot{X}^T_{T/2}+F'(X_{T/2}^T),\dot{X}^T_{T/2}-F'(X_{T/2}^T) \rangle| \leq |\dot{X}^T_{T/2}+F'(X_{T/2}^T)||\dot{X}^T_{T/2}-F'(X_{T/2}^T)|.
		\end{equation*}
		It follows from Proposition \ref{ggprop-2} that
		
		\begin{equation*}
		|\dot{X}^T_{T/2}+F'(X_{T/2}^T)|\leq \sqrt{\frac{2\rho}{\exp(\rho T)-1}(\mathcal{C}_T(x,y)+2F(y)-2F(x) ).}
		\end{equation*}
		Next, we observe that the time-reversal of $(X^T_t)_{t\in[0,T]}$ is optimal for the variational problem obtained exchanging the labels $x$ and $y$ in \eqref{eq-48}. This implies that $C_{T}(x,y)=C_{T}(y,x)$ and thanks again to Proposition \ref{ggprop-2} that
		
		\begin{equation*}
		|\dot{X}^T_{T/2}-F'(X_{T/2}^T)|\leq \sqrt{\frac{2\rho}{\exp(\rho T)-1}(\mathcal{C}_T(x,y)+2F(x)-2F(y) ).}
		\end{equation*}
		using these two bounds in the above expression gives \eqref{eq-202}.
	\end{eproof}
	
	We are now ready to prove the announced result. The proof is based on the above proposition and the finite-dimensional version of the logarithmic Sobolev inequality, which reads as
	
	\begin{equation}\label{eq-203}
	2\rho ( F(x)-F(x^*)) \leq |F'(x)|^2, \quad \forall x\in \R^n.
	\end{equation}
	
	\begin{ethm} \label{gio-thm2}
		For all $x,y>0$, $T>0$ and $t\in(0,T)$ we have:

		\begin{multline}
		\label{eq-204}
		F(X_t^T) \leq \frac{\sinh(2\rho(T-t))}{\sinh(2\rho T)}\left(F(x)-\frac{E_T(x,y)}{4 \rho}+F(x^*) \right) \\ +\frac{\sinh(2\rho t)}{\sinh(2\rho(T-t))}\left(F(y) - \frac{E_T(x,y)}{4 \rho}+F(x^*) \right)+\frac{E_T(x,y)}{4 \rho}-F(x^*)
		\end{multline}
		Moreover, for all fixed $\theta\in(0,1)$ there exists a decreasing function $b(\cdot)$ such that 
		
		\begin{equation}\label{eq-205}
		F(X^{T}_{\theta T})-F(x^*) \leq b(\rho)(F(x)+F(y)-2F(x^*)) \exp(-2\rho \min\{\theta,1-\theta\}T). 
		\end{equation}
		holds uniformly in $T\geq 1$.

	\end{ethm}
	
	\begin{eproof}
		A standard calculation gives
		
		\begin{equation*}
		\frac{d}{dt}F(X^{T}_t) = \langle F'(X_t^T),\dot{X}_t^T\rangle =  \frac{1}{4}\big( |F'(X_t^T)+\dot{X}_t^T|^2- |F'(X_t^T)-\dot{X}_t^T|^2 \big)
		\end{equation*}
		
		From this expression we obtain, using Newton's equation and $\rho$-convexity of $F$:
		
		\begin{align*}
		\frac{d}{dt}\frac{1}{4}\big( |F'(X_t^T)+\dot{X}_t^T|^2- |F'(X_t^T)-\dot{X}_t^T|^2 \big) \geq \frac{\rho}{2} (|F'(X_t^T)+\dot{X}_t^T|^2+|F'(X_t^T)-\dot{X}_t^T|^2 )\\
		=2\rho|F'(X^T_t)|^2 +\rho E_{T}(x,y).
		\end{align*}
		
		At this stage we can use the logarithmic Sobolev inequality~\eqref{eq-203} to obtain that 
		$$
		2\rho|F'(X^T_t)|^2 +\rho E_{T}(x,y) \geq 4\rho^2 \left( F(X^{T}_t) - F(x^*) \right)+\rho E_{T}(x,y).
		$$
		Summing up, we have obtained that the function $t\mapsto F(X^{T}_t)$ satisfies the differential inequality
		$$
		\frac{d^2}{dt^2}F(X^{T}_t) \geq 4\rho^2 \left( F(X^{T}_t) -F(x^*)\right) +\rho E_{T}(x,y).
		$$
		The bound \eqref{eq-205} is then obtained integrating this differential inequality, see \cite[Lemma 5.6]{backhoff2019mean} for details. The bound \eqref{eq-205} follows by using \eqref{eq-202} and the upper bound \eqref{eq-206} in \eqref{eq-204} after some standard (though tedious) calculations.
	\end{eproof}

	\subsection{The $(0,n)$-convex case}
	\label{sec-3.3}

	Now we assume an other kind of convexity. We assume $F$ is $(0,n)$-convex that is
	\begin{equation}
	\label{eq-400}
	F'' \geq \frac{1}{n}F' \otimes F'.
	\end{equation}

	\subsubsection{Costa type estimates under the $(0,n)$-convexity}
	\label{sec-3.3.1}
	The $(0,n)$-convexity is related to Costa type convexity~\cite{costa} and produced many useful estimates.  All estimates  are related to the same trick. 
	Let  $a > 0$, $T>0$ and $\phi:[0,T] \rightarrow \R$ a smooth function satisfying
	\begin{equation}
	\label{eq-77}
	\forall t\in[0,T],\quad \frac{d}{dt}\phi(t)\geq a \phi^2(t).
	\end{equation}
	Let $\Phi$ be an antiderivative of $\phi$, then the map $\Lambda(t)= e^{-a \Phi(t)}$, ($t\in[0,T]$) is a concave function on $[0,T]$. In that case, coming from classical convex inequalities for the function $\Lambda$, 
	$$
	\Lambda'(T)\leq\frac{\Lambda(T)-\Lambda(t)}{T-t}\leq \Lambda'(t)\leq\frac{\Lambda(t)-\Lambda(0)}{t}\leq \Lambda'(0),
	$$
	one can deduce the following properties.
	\begin{enumerate}[1.]
		\item For all $t\in(0,T)$,
		\begin{equation}
		\label{eq-64}
		-\frac{1}{at}\leq \phi(t)\leq \frac{1}{a(T-t)}.
		\end{equation}
		\item We also have the following inequality
		\begin{equation}
		\label{eq-63}
		-\frac1a\log (1-aT\phi(0))\leq\Phi(T)-\Phi(0)\leq \frac1a\log (1+aT\phi(T)).
		\end{equation}
	\end{enumerate}
	
	In our case, this remark gives some important estimates for gradient flow or $F$-interpolation where the proofs are elementary. 
	\begin{enumerate}
		\item Costa's convexity~\cite{costa}: for any $x\in\R^n$, the map
		\begin{equation}
		\label{eq-60}
		[0,\infty)\ni t\mapsto \exp\PAR{-\frac{2}{n}F(S_t(x))}
		\end{equation}
		is concave. Recall that $(S_t(x))_{t \geq 0}$ is the gradient flow of $F$ with initial position $x$, defined in~(\ref{eq-1001}). 
		\item Ripani's convexity~\cite{ripani}: for any $F$-interpolation $(X_t^T)_{t\in[0,T]}$, the map 
		\begin{equation}
		\label{eq-61}
		[0,\infty)\ni t\mapsto \exp\PAR{-\frac{1}{n}F(X_t^T)}
		\end{equation}
		is concave. 
		\item Improved Ripani's convexity: for any $F$-interpolation $(X_t^T)_{t\in[0,T]}$, the map
		\begin{equation}
		\label{eq-62}
		[0,\infty)\ni t\mapsto \exp\PAR{-\frac{1}{n}\SBRA{F(X_t^T)+\int_0^t|F'(X_s^T)|^2ds}}
		\end{equation}
		is concave.
	\end{enumerate}

	\subsubsection{Convergence of the $F$-interpolation}
	\label{sec-3.3.2}
	
	\black

	{
		
		We begin by proving that the derivative of the cost in $T$ is precisely $- E_{T}(x,y)$, as observed in~\cite{conforti-tamanini2019} for the classical Schrödinger problem.
		\begin{eprop}
			We have for all $x,y\in\R^n$ and $T>0$ that
			
			\begin{equation}\label{eq-207}
			\frac{d }{d T} C_{T}(x,y) = - E_{T}(x,y).
			\end{equation}
		\end{eprop}
		
		\begin{eproof}
			Here we need to introduce another formulation of the cost. For every $x,y \in \R^n$ and $T>0$ we define
			$$
			A_T(x,y)=\inf \int_0^1 \big[|\dot \omega_s|^2+T^2|F'(\omega_s)|^2\big]ds,
			$$
			where the infimum runs over all paths from $x$ to $y$. Then from the so called envelope theorem (see e.g. \cite{lafrance1991envelope} for a formulation of the envelope in the context of dynamic control problems) and recalling that $A_{T}(x,y)=TC_T(x,y)$ we obtain
			\begin{equation*}
			\frac{d}{dT}A_{T}(x,y)=\frac{d}{dT}TC_{T}(x,y) = 2 T\int_{0}^{1}|F'(\tilde{\omega}_s)|^2ds,
			\end{equation*}
			where $\tilde{\omega}$ is the optimal path in $A_{T}(x,y)$. Operating the change of variable $Tt=s$ we get that 
			\[ 
			\frac{d}{dT}TC_{T}(x,y) = 2\int_{0}^{T}|F'(X^T_t)|^2 dt,
			\]
			where $X^T_t$ is the $F$ interpolation between $x$ and $y$. Adding and substracting $|\dot{X}^T_t|^2$ in the integral and observing that the definition of cost and conserved quantity we arrive at 
			
			\begin{equation*}
			\frac{d}{dT}TC_{T}(x,y)=C_{T}(x,y)-TE_{T}(x,y),
			\end{equation*}
			from which the desired conclusion follows.
		\end{eproof}
		
		As in the $\rho$-convex case we introduce $\varphi_t^T= \dot X_t^T +F'(X_t^T)$. Combining the latter with the improved Ripani convexity yields some useful results.

		\begin{ethm}
			For any $x,y\in\R^n$ and $T>0$ we have 
			\begin{equation}\label{eq-208}
			-E_{T}(x,y) \leq \frac{2n}{T},\quad C_{T}(x,y) \leq C_{1}(x,y) + 2n \log T,
			\end{equation}
			Moreover, for all $t\in(0,T)$ we have
			\begin{equation}
			|\varphi^T_t|^2 \leq \frac{2F(y)-2F(x) +C_{1}(x,y) + 2n \log T}{T-t},
			\end{equation}
			where $\varphi_t^T:= \dot X_t^T+F'(X_t^T)$ for every $T>0$ and $t \in [0,T]$ .
			
		\end{ethm}
		
		\begin{eproof}
			For the first statement observe that by \eqref{eq-62} and the trick we have 
			$$
			|F'(X_0^T)|^2+\langle \dot X_0^T,F'(X_0^T) \rangle \leq \frac {n}{T},
			$$  
			and completing the squares we have that
			\begin{equation*}
			|F'(X_{0}^T)|^2+\langle F'(X_{0}^T),\dot X_{0}^T\rangle \geq \frac{1}{2}|F'(X_{0}^T)|^2-\frac{1}{2}|\dot X_{0}^T|^2=-\frac12E_{T}(x,y),
			\end{equation*} 
			which gives the first bound $-E_T(x,y) \leq \frac{2n}{T}$.  Integrating this inequality between $1$ and $T$ we find the desired bound for the cost. Finally, for the last inequality  we observe that, $\frac{d}{dt} |\varphi_t^T|^2=F''(X_t^T)(\varphi_t^T,\varphi_t^T)$, hence the function $t \mapsto |\varphi_t^T|^2$ is non decreasing, hence 
			
			\begin{align*} (T-t)|\varphi_t^T|^{2} &\leq \int_{t}^{T}|\varphi_s^T|^{2} d s\leq \int_{0}^{T}|\varphi_s^T|^{2} d s \\
			&= \int_{0}^{T}|\dot{X}^T_s|^2 +2\langle F'(X^T_s),\dot{X}^T_s\rangle + |F'(X^T_s)|^2 d s \\
			&=C_{T}(x,y)+ 2F(y)-2F(x).
			\end{align*}
			Using \eqref{eq-208} we get the desired result.
		\end{eproof}
		\newline
		
		Note that since $T|E_{T}(x,y)|  \leq C_{T}(x,y)$ we also obtain the two-sided bound $|E_{T}(x,y)|\leq \log(T)/T$. As in the $\rho$-convex case, we can deduce the convergence of the $F$-interpolation towards the gradient flow semigroup of $F$ from the estimate of $|\phi_t^T|^2$. The proof is exactly the same as in the $\rho$-convex case, using the previous estimate.
		
		\begin{ethm}[Distance between entropic interpolations and gradient flows]
			For all $x,y \in \R^n$, $T>2$ and $t \in (0,T)$ we have $$
			|X_t^T-S_t(x)| \leq 2\sqrt{\left(2(F(y)-F(x))+C_1(x,y)+2n \log(T)\right)}\left(\sqrt{T}- \sqrt{T-t} \right).
			$$
			In other words,  for all $a>2$, there exists a constant $C \geq 0$ such that for all $T > a$ and $t \in [0,a]$
			\begin{equation} \label{eq-65}
			|X_t^T-S_t(x)| \leq C  \sqrt{\frac{n \log(T)}{T}}.
			\end{equation}
		\end{ethm}
		In light of the example described in Section~\ref{sec-3.2.1}, the estimate~\eqref{eq-65} may not be optimal.
		
		\subsubsection{A turnpike estimate}
		We saw in section \ref{sec-3.2} that the fundamental exponential entropy dissipation estimate along the heat flow can be generalized to $F$-interpolations under the $CD(\rho,\infty)$ condition. Under the condition $CD(0,n)$ the following fundamental estimates for the Fisher information along the heat flow is known to hold,
		\begin{equation*}
		|F'(S_t(x))|^2 \leq \frac{n}{2t}.
		\end{equation*}
		To prove such inequality, is it enough to differentiate $|F'(S_t(x))|^2$ in time and apply the $(0,n)$-convexity property of $F$ to close a differential inequality. 
		In the next result  we generalize this estimate to $F$-interpolations. It is worth noticing that Theorem 3.8 below  yields meaningful information at timescales that are $O(1)$, i.e. when $t$ is fixed. On the contrary, the next result yields a non trivial bound also at timescales that are of the order $O(T)$. 
		
		\begin{ethm}\label{gio-thm3}
			For any $x,y\in\R^n$, $T>0$ and $t\in(0,T)$ we have 
			
			\begin{equation}\label{e-210}
			|F'(X^T_{t})|^2 \leq   \frac{n}{{{2}}t}+\frac{n}{{{2}}(T-t)},
			\end{equation}
			furthermore for every $T>0$ and $\theta \in (0,1)$,
			$$
			|F'(X_{\theta T}^T)|^2 \leq \frac{n}{{{2}}T \theta (1- \theta)}.
			$$
		\end{ethm}
		
		\begin{eproof}
			The proof consists in combining inequalities of Section~\ref{sec-3.3.1} and time-reversal. We first observe that from improved Ripani convexity~\eqref{eq-62} we have
			\begin{equation}\label{e-170}|F'(X^T_{t})|^2+\langle F'(X_t^T),  \dot X_t^T \rangle\leq \frac{n}{T-t}.\end{equation}
			Next, we remark that $(Y^{T}_t)_{t\in[0,T]}=(X^{T}_{T-t})_{t\in[0,T]}$ is a $F$-interpolation between $y$ and $x$, i.e. it is optimal for the variational problem obtained from \eqref{eq-48} inverting the roles of $x$ and $y$. But then, using again~\eqref{eq-62},
			$$
			|F'(Y^T_{T-t})|^2+\langle F'(Y_{T-t}^T),  \dot Y_{T-t}^T \rangle\leq \frac{n}{t}.
			$$
			which is equivalent to
			$$ |F'(X^T_{t})|^2-\langle F'(X_{t}^T),  \dot X_{t}^T \rangle\leq \frac{n}{t}.$$
			Adding up this last bound and \eqref{e-170} yields the desired result.
		\end{eproof}
		
		
		

		\section{The infinite dimensional case}\label{sec:main}
		\label{sec-4} 
		
		From now on, our base space the space of probability measure $\mathcal P_2(N)$ instead of $\R^n$. In what follows, we shall see how it is possible to replicate in a rigorous fashion the results obtained in the finite dimensional case in the infinite dimensional setup. 
		
		\subsection{The example of two Gaussian measures on $\R$}
		
		As we did before, we perform some explicit calculation, using some simple example in order to build intuition. In this example $N=\R$ is the Euclidean space equipped with the classical Laplace operator.  We are gonna to compute all the desired quantities in the case of two Gaussian measures in $\R$. Let $x_0,x_1 \in \R$, $\mu= \mathcal{N}(x_0,1)$ and $\nu = \mathcal{N}(x_1,1)$. We denote by $\mathcal{N}(m,\sigma^2)$ the usual Gaussian distribution with mean $m$ and variance $\sigma^2$.  
		
		Recall that the gradient flow of the standard entropy is the dual of the classical heat semigroup, namely
		$$
		\mathcal S_t (\mu)=P_t^{*}(\mu)=P_t \left(\frac{d\mu}{d Vol} \right)d Vol, \ t \geq 0.
		$$
		In this particular case we are able to compute all the quantities of interest.
		\begin{itemize}
			\item The gradient flow starting from $\mu$ of the standard entropy is given by
			$$
			\forall t>0,\, \mathcal S_t(\mu)=P_t\left(\frac{d\mu}{d Vol}\right)d Vol=\mathcal{N}\left(x_0,1+2t \right).
			$$
			\item The entropic interpolation between $\mu$ and $\nu$ is the path $(\mathcal{N}\left(x_t^T, \sigma_t^T\right))_{t\in[0,T]}$ where
			\begin{equation*}
			\left\{\begin{array}{l}
			\displaystyle  x_t^T=\frac{T-t}{T}x_0+\frac{t}{T}x_1,  \\
			\displaystyle \sigma_t^T=1+2\frac{t(T-t)}{D_T^2+T},
			\end{array} \right.
			\end{equation*}
			for some $D_T>0$. Furthermore it can be shown that, for every $t>0$,
			\begin{equation*}
			\left\{\begin{array}{l}
			\displaystyle  x_t^T \underset{T \rightarrow \infty}{\rightarrow} x_0,  \\
			\displaystyle \sigma_t^T \underset{T \rightarrow \infty}{\rightarrow} 1+2t.
			\end{array} \right.
			\end{equation*}

			\item  The conserved quantity is given by
			$
			\mathcal E_T(\mu,\nu)=\frac{1}{T^2}(x_1-x_0)^2- \frac{2}{(\mathcal{D}_T^2+T)},
			$
			that is, for some constant $c>0$,
			$$
			\left|\mathcal E_T(\mu,\nu) \right| \leq \frac{c}{T}.
			$$
			\item  We have the following estimate for the cost 
			$
			\mathcal{C}_T(\mu,\nu) \underset{T \rightarrow \infty}{\sim} 2 \log (T).
			$
			
			\item The distance between the entropic interpolations and the gradient flow is given by, with some $C>0$,
			$$
			W_2\left(\mu_t^T,P_t^{*} \mu\right)=\sqrt{\left|\sqrt{1+2t}-\sqrt{\sigma_t^T}\right|^2+\left|x_0-x_t^T\right|^2 }\underset{T \rightarrow \infty}{\sim}  \frac{C}{T}.
			$$
		\end{itemize}
		As a conclusion, at least for two Gaussian measures, we obtain the same asymptotic behaviour as in the example given in Section~\ref{sec-3.2.1}.

		\subsection{The $CD(\rho,\infty)$ case}
		
		In this subsection we assume that the semigroup $(P_t)_{t \geq 0}$ verifies a $CD(\rho, \infty)$ curvature-dimension condition that is
		\begin{equation*}
		\Gamma_2(f)\geq \rho\Gamma(f),
		\end{equation*}
		with $\rho>0$. Since $\rho>0$, the reversible measure $m$ is then a probability measure and the functional $\mathcal F $ has a unique minimum $m$ such that $\mathcal F(m)=0 $.

		First, as in the finite dimensional case, a Talagrand type inequality for the entropic cost who gave a bound for the cost, that is for all $\mu,\nu$ compactly supported probability measures,
		$$
		\mathcal C_T(\mu,\nu) \leq 2 \inf_{t\in(0,T)}\left\{ \frac{1+e^{-2 \rho t}}{1-e^{-2\rho t}}\mathcal{F}(\mu)+\frac{1+e^{-2\rho(T-t)}}{1-e^{-2\rho(T-t)}} \mathcal{F}(\nu)\right\}.
		$$
		In particular we have the following inequality
		$$
		\mathcal{C}_T(\mu,\nu) \leq 2 \frac{1+e^{- \rho T}}{1-e^{- \rho T}}\left(\mathcal{F}(\mu) + \mathcal{F}(\nu) \right).
		$$

		These inequalities were first obtained in~\cite{conforti2017} (see also~\cite[Cor~4.5]{gentil-leonard2020}). Hence the entropic cost is bounded under a $CD(\rho,\infty)$ condition. As in the finite dimensional case we need an estimate for the $L^2(\mu_t^T)$ norm of
		$$
		\varphi_t^T:= \grad_{\mu_t^T} \mathcal{F}+ \dot \mu_t^T= \nabla \log \left(\frac{d \mu_t^T}{dm} \right)+\dot \mu_t^T, \, t \in [0,T].
		$$
		The following proposition, whose proof has been implicitly already done in~\cite[Thm 1.4]{conforti-tamanini2019} is  obtained following the proof of Proposition~\ref{ggprop-2}. 
		
		\begin{eprop} 
			\label{gioprop-1} Let assume the $CD(\rho,\infty)$ condition with $\rho>0$. 
			Let $\mu,\nu \in \mathcal{P}_2(N)$ be two absolutely continuous, compactly supported measures with smooth positive densities against $m$. For all $T>0$, $\left( \mu_t^T\right)_{t \in [0,T]}$ denotes the entropic interpolation from $\mu$ to $\nu$ and for $t \in (0,T)$ we define $\varphi_t^T:= \dot \mu_t^T+ \nabla \log (\mu_t^T)$. Then for all $T>0$ and $t \in (0,T)$
			$$
			|\varphi_t^T|_{\mu_t^T}^2 \leq  \frac{2 \rho}{\exp(2\rho(T-t))-1}\left(\mathcal{C}_T(\mu,\nu)+2\mathcal{F}(\nu)-2\mathcal{F}(\mu) \right).
			$$
			In particular
			$$
			|\varphi_t^T|_{\mu_t^T}^2 \leq \frac{8 \rho}{\exp(2 \rho  (T-t))-1}\left(\frac{e^{-2 \rho t}}{1-e^{-2 \rho t}}\mathcal{F}(\mu)+\frac{1}{1-e^{-2 \rho (T-t)}}\mathcal{F}(\nu) \right).
			$$
		\end{eprop}
		
		Now we can obtain the main result: convergence of entropic interpolation towards gradient flow. The idea of the proof is the same as in the finite dimensional case, however, some extra care has to be taken in order to differentiate the Wasserstein distance.  We recall that if $(\delta_t)_{t\geq 0}$ and $(\eta_t)_{t\geq 0}$ be two absolutely continuous curves in $\mathcal P_2(N)$ such that for every $t\geq 0$, $\delta_t$ and $\eta_t$ are absolutely continuous w.r.t. $d Vol$.  Then we have for almost every $t\geq 0$,
		$$
		\frac{d}{dt} \frac{ W_2^2(\delta_t, \eta_t)}{2}=-\langle T^1_t,\dot \delta_t \rangle_{\delta_t}-\langle  T^2_t, \dot \eta_t \rangle_{\eta_t},
		$$
		where $\exp(T^1_t)$ (resp. $\exp(T^2_t)$) is the optimal transport $\delta_t \rightarrow  \eta_t$ (resp. $\eta_t \rightarrow  \delta_t$), see~\cite[Theorem~23.9]{villani2009}. Now we can state our main theorem. 
		
		\begin{ethm}[Convergence of the entropic interpolation]
			\label{gthm-2}
			Let assume the $CD(\rho,\infty)$ condition with $\rho>0$. Let $\mu,\nu \in \mathcal P_2(N)$ be two absolutely continuous, compactly supported measures with smooth positive densities w.r.t.  $m$. For all $T>0$, $\left( \mu_t^T\right)_{t \in [0,T]}$ denotes the entropic interpolation from $\mu$ to $\nu$.
			Then for all $T>0$ and $t \in (0,T)$,
			$$
			W_2\left(\mu_t^T,P_t^{*}(\mu) \right) \leq t\exp(- \rho T ) \sqrt{\frac{2 \rho}{\exp(-2\rho t)-\exp(-2\rho T)}\left(\mathcal{C}_T(\mu,\nu)+2(\mathcal{F}(\nu)-\mathcal{F}(\mu)) \right)}.
			$$ 
			In other words, there exists a constant $C>0$ depending only on $\mu$ and $\nu$ such that for every $t\geq 0$ and $T>t$,
			$$
			W_2(\mu_t^T,P_t^{*}\mu) \leq C \frac{t \exp(-\rho T)}{\sqrt{\exp(-2\rho t)-\exp(-2 \rho T)}}.
			$$
			
		\end{ethm}

		\begin{eproof}
			Let $T>0$ and $t$ be a Lebesgue point of $[0,T]$. The derivative of the  Wasserstein distance gives
			\begin{equation*}
			\begin{split}
			\frac{d}{dt}\frac{ W_2^2(\mu_t^T,P_t^*\mu)}{2} &=-\langle T^1_t, \dot \mu_t^T \rangle_{\mu_t^T}+\langle T^2_t,  \grad_{P_t^*\mu} \mathcal F\rangle_{P_t^*\mu} \\
			&= -\langle T^1_t , \dot \mu_t^T+ \grad_{\mu_t^T} \mathcal F \rangle_{\mu_t^T}+ \langle T^1_t, \grad_{\mu_t^T} \mathcal F \rangle_{\mu_t^T}+ \langle T^2_t, \grad_{P_t^* \mu} \mathcal{F} \rangle_{P_t^* \mu} .
			\end{split}
			\end{equation*}
			Where $\exp(T^1_t)$ (resp. $\exp(T^2_t)$) is the optimal transport $\mu_t^T\rightarrow P_t^* \mu$ (resp.  $P_t^*\mu\rightarrow\mu_t^T$). From~\cite[Theorem~23.14]{villani2009} we have 
			$$
			\langle T^1_t, \grad_{\mu_t^T} \mathcal F \rangle_{\mu_t^T}+ \langle T^2_t, \grad_{P_t^* \mu} \mathcal{F} \rangle_{P_t^* \mu} \leq 0,
			$$
			which is actually a rigorous proof of the convexity of the entropy along a geodesic. Whence we have obtained  
			$$
			\frac{d}{dt}\frac{ W_2^2(\mu_t^T,P_t^* \mu)}{2} \leq  - \langle T^1_t,\dot \mu_t^T + \grad_{P_t^* \mu} \mathcal F \rangle_{\mu_t^T}. 
			$$
			
			Since $|T^1_t|_{\mu_t^T}=W_2(\mu_t^T,P_t^* \mu)$, by the Cauchy-Schwarz inequality
			$$
			\frac{d}{dt} \frac{ W_2^2(\mu_t^T,P_t^*\mu)}{2} \leq |\varphi_t^T|_{\mu_t^T}W_2(\mu_t^T,P_t^*\mu),
			$$
			where $\varphi_t^T=\dot\mu_t^T+\grad_{P_t^* \mu} \mathcal{F}$. These inequalities holds for almost every $t\in [0,T]$. The result follow from the integration of this inequality and the Proposition~\ref{gioprop-1}, as in the finite dimensional case.
		\end{eproof}

		\subsubsection{Turnpike property}
		
		It is well known that under the $CD(\rho,\infty)$ curvature dimension condition the gradient flow $P_t^*$ of $\mathcal{F}$ dissipates at exponential rate $2\rho$, in particular for $T>0$ and $\theta \in (0,1)$ we have for $\mu \in \mathcal{P}_2(N)$
		$$
		\mathcal{F}(P_T^*(\mu)) \leq \mathcal{F}(\mu)\exp(-2\rho T).
		$$
		Recall that in this case $\mathcal F(m)=0$ and is the minimum of $\mathcal F$ on $\mathcal P(N)$. 
		
		As in the finite dimensional case we can show that a similar estimate holds along entropic interpolations. The first step is an exponential upper bound for the conserved quantity. The proof is exactly the same as in the finite dimensional case. 
		
		\begin{eprop} Let $\mu,\nu \in \mathcal{P}_2(N)$ be two compactly supported absolutely continuous measures with smooth positive densities w.r.t. $m$. Then for every $T>0$
			$$
			|\mathcal{E}_T(\mu,\nu)| \leq \frac{2 \rho}{\exp(\rho T)-1}\sqrt{\mathcal{C}_T^2(\mu,\nu)-4(\mathcal{F}(\mu)-\mathcal{F}(\nu))^2}.
			$$
		\end{eprop}
		
		We can now state our main result, the proof is similar to the proof of the Theorem~\ref{gio-thm2} using the logarithmic Sobolev inequality. A similar result has been obtained for the mean field Schr\"odinger problem, see \cite{backhoff2019mean}.
		$$
		2\rho \mathcal{F}(\mu) \leq |\grad_{\mu} \mathcal{F}|^2_{\mu}=\mathcal{I}_W(\mu).
		$$
		
		\begin{ethm} Let $\mu,\nu \in \mathcal{P}_2(N)$ be two be two compactly supported absolutely continuous measures with smooth positive densities w.r.t.  $m$. Then For every $T>0$ and $t \in (0,T)$ we have
			$$
			\mathcal{F}(\mu_t^T) \leq \frac{\sinh(2 \rho (T-t))}{\sinh(2\rho T)}\left(\mathcal{F}(\mu)-\frac{\mathcal{E}_T(\mu,\nu)}{4\rho}\right)+\frac{\sinh(2 \rho t)}{\sinh(2\rho T )}\left(\mathcal{F}(\nu)-\frac{\mathcal{E}_T(\mu,\nu)}{4\rho}\right) +\frac{\mathcal{E}_T(\mu,\nu)}{4\rho}.
			$$
			Moreover, for all $\theta \in (0,1)$ there exists a decreasing function $b(\cdot)$ such that
			$$
			\mathcal{F}(X_{\theta T}^T) \leq b(\rho)(\mathcal{F}(\mu)+\mathcal{F}(\nu))\exp(-2 \rho \min \left\{\theta, 1 - \theta\right\}T).
			$$
			
		\end{ethm}

		\subsection{The $CD(0,n)$ case}
		\label{sec-4.3}
		
		In this subsection we assume that $L$ satisfies a $CD(0,n)$ curvature-dimension condition, that is for every smooth function $f$,
		$$
		\Gamma_2(f) \geq \frac{1}{n} \left(L P_tf \right)^2.
		$$
		As explained in Section~\ref{sec:setting}, this case is covers the fundamental example of $\R^n$ with the usual  Laplacian. In that case, the measure $m$ is not a probability measure.

		The aim of this subsection is to prove the convergence of entropic interpolations towards the semigroup $(P_t^*)_{t \geq 0}$ under the $CD(0,n)$ condition. But first let's recall Costa type estimates, which are fundamental for our purpose. Some of these results are generalized in~\cite{ccg2}.
		
		\subsubsection{Costa type estimates under the $CD(0,n)$ condition}
		\label{sec-4.3.1}
		
		The $CD(0,n)$ condition gives some important estimates for gradient flow or entropic interpolations. The proofs follow the same trick explained in Section~\ref{sec-3.3.1}. As in the finite dimensional case,   estimates are given for the gradient flow or the entropic interpolation. 
		\begin{enumerate} 
			\item Costa's convexity~\cite{costa} : for any $\mu\in \mathcal P_2(N)$ the map 
			\begin{equation}
			\label{eq-70}
			t\ni[0,\infty)\mapsto \exp\PAR{-\frac{2}{n}\mathcal F(P_t^*(\mu))}
			\end{equation}
			is concave.
			
			Let us briefly  recall the proof. For any probability measure $\mu$, 
			$$
			\mathcal F(P_t^*(\mu))=\int P_t h\log P_t hdm=\mathcal Ent_m(P_th),
			$$
			where $h=\frac{d\mu}{dm}$. Following the Bakry-\'Emery computations, see for instance~\cite[Proof of Theorem~6.7.3]{bgl-book} 
			\begin{multline*}
			\frac{d^2}{dt^2}\mathcal F(P_t^*(\mu))=2\int \Gamma_2\PAR{\log P_th}P_thdm\geq \frac{2}{n}	\int \PAR{L\log P_th}^2P_thdm\geq\\\frac{2}{n}	\PAR{\int {L\log P_th}P_thdm}^2= \frac{2}{n}	\Big(\int \Gamma\PAR{\log P_th}P_thdm\Big)^2=
			\frac{2}{n}	\Big(\frac{d}{dt}\mathcal F(P_t^*(\mu))\Big)^2,
			\end{multline*}
			which is the inequality~\eqref{eq-77} with $a=2/n$.
			
			As in the finite dimensional case we obtain two inequalities useful for the rest of the paper,
			\begin{equation}
			\label{eq-80}
			\mathcal I_W(P_t^*\mu)\leq \frac{2}{nt},
			\end{equation}
			and 
			\begin{equation}
			\label{eq-81}
			\mathcal F(\mu)-\mathcal F(P_T^*(\mu))\leq\frac n2\log \PAR{1+\frac{2T}{n}\mathcal I_W(\mu)}.
			\end{equation}

			\item Ripani's convexity~\cite{ripani}: for any entropic interpolation $(\mu_t^T)_{t\in[0,T]}$, the map
			\begin{equation}
			\label{eq-71}
			t\ni[0,\infty)\mapsto \exp\PAR{-\frac{1}{n}\mathcal F(\mu_t^T)}
			\end{equation}
			is concave.
			
			The proof is similar to Costa's convexity. There exits two positive functions $f,g$ such that
			$$
			\mu_t^T=P_tfP_{T-t}g \,m, 
			$$ 
			then 
			$$
			\mathcal F(\mu_t^T)=\int P_t fP_{T-t}g\log (P_t fP_{T-t}g)dm, 
			$$
			and the proof is based on computation of the second derivative of such function, see~\cite{ripani} for additional details.
			
			\item Improved Ripani's convexity: for any entropic-interpolation $(\mu_t^T)_{t\in[0,T]}$, with $\mu_0^T$ and $\mu_T^T$ smooth and compactly supported probability measures, the map
			\begin{equation}
			\label{eq-72}
			t\ni[0,\infty)\mapsto \exp\PAR{-\frac{1}{n}\SBRA{\mathcal F(\mu_t^T)+\int_0^t|{\rm grad}_{\mu_s^T}\mathcal F|_{\mu_s^T}^2ds}}
			\end{equation}
			is concave. 
			
			In particular, from~\eqref{eq-64}, we obtain for $t\in[0,T)$, 
			\begin{equation}
			\label{eq-82}
			\langle\grad_{\mu_t^T}\mathcal F,\dot\mu_t^T\rangle+|{\rm grad}_{\mu_t^T}\mathcal F|_{\mu_t^T}^2=   \int 2P_t\PAR{\frac{\Gamma(P_{T-t}g)}{P_{T-t}g}-L P_{T-t}g}f\,dm\leq \frac{n}{T-t}
			\end{equation}

			A rigorous proof of the concavity of~\eqref{eq-72} is quite tricky. For an heuristic proof, it is enough to formally compute the second derivative of~\eqref{eq-72} and use the infinite dimensional version of~\eqref{eq-400}. It will be discussed in a forthcoming paper. For the scope of this paper, we only need inequality~\eqref{eq-82} for which we can provide a direct proof.  Again, there exist two positive smooth and compactly supported functions $f,g$ such that $\mu_t^T=P_tfP_{T-t}g \,m$, then 
			$$
			\grad_{\mu_t^T}\mathcal F=\nabla \log (P_tfP_{T-t}g),
			$$
			and 
			$$
			{\dot\mu_t^T}=\nabla \log P_{T-t}g-\nabla \log P_tf. 
			$$
			Then we obtain,
			\begin{multline*}
			\langle\grad_{\mu_t^T}\mathcal F,\dot\mu_t^T\rangle+|{\rm grad}_{\mu_t^T}\mathcal F|_{\mu_t^T}^2=\\
			\int \Big(\Gamma\Big(\log(P_tfP_{T-t}g),\log\frac{P_{T-t}g}{P_tf}\Big)+\Gamma(\log (P_tfP_{T-t}g))\Big)P_tfP_{T-t}g\,dm=\\
			\int 2P_t\PAR{\frac{\Gamma(P_{T-t}g)}{P_{T-t}g}-L P_{T-t}g}f\,dm.
			\end{multline*}
			The so-called Li-Yau inequality, proved for instance in~\cite{bakry-ledoux} in the context of the $CD(0,n)$-condition, insures that for $t\in[0,T)$,
			$$
			\frac{\Gamma(P_{T-t}g)}{(P_{T-t}g)^2}-\frac{L P_{T-t}g}{P_{T-t}g}\leq \frac{n}{2(T-t)},
			$$
			which implies~\eqref{eq-82}.

		\end{enumerate}
		
		
		\subsubsection{Convergence of the entropic interpolation}
		
		In this subsection we follow exactly the line of reasoning adopted in the finite dimensional $(0,n)$-convex case. We first notice that the derivative of the cost in $T$ is exactly $-\mathcal{E}_T(\mu,\nu)$.
		
		\begin{eprop}\label{prop: envelope}[\cite{conforti-tamanini2019}] Let $\mu,\nu \in \mathcal{P}_2(N)$  be two absolutely continuous and compactly supported measures with smooth density w.r.t. $m$. Then for every $T>0$
			$$
			\frac{d}{dT}\mathcal{C}_T(\mu,\nu)=-\mathcal{E}_T(\mu,\nu).
			$$
		\end{eprop}
		
		Defining $\varphi_t^T:=\dot \mu_t^T+ \grad_{\mu_t^T}\mathcal{F}= \dot \mu_t^T + \nabla \log \left(\frac{d \mu_t^T}{dm} \right)$, combining the latter with Ripani convexity we obtain, exactly as in the finite dimensional case, the following result.
		
		\begin{ethm}[Large time asymptotics for cost and energy]\label{thm-cost-consqty-asym} Let $\mu,\nu \in \mathcal{P}_2(N)$ be two  compactly supported absolutely continuous measures with smooth positive densities w.r.t. $m$. For all $T>1$, we denote by $\left( \mu_t^T\right)_{t \in (0,T)}$ the entropic interpolation from $\mu$ to $\nu$ and for $t \in (0,T)$ we define $\varphi_t^T:= \dot \mu_t^T+ \nabla \log (\mu_t^T)$. Then for every $T>0$ we have
			$$
			-\mathcal{E}_T(\mu,\nu) \leq \frac{2n}{T}, \quad \mathcal{C}_T(\mu,\nu) \leq \mathcal{C}_1(\mu,\nu)+2n\log(T),
			$$
			and for all $t \in (0,T)$ we have
			\begin{equation}\label{eq-305}
			|\varphi_t^T|_{\mu_t^T}^2 \leq \frac{2 \mathcal{F}(\nu)-2\mathcal{F}(\mu)+\mathcal{C}_1(\mu,\nu)+2n \log(T)}{T-t}.
			\end{equation}
		\end{ethm}

		Now we can state the main result of this subsection. The proof is the exact analogous of Proposition~\ref{gthm-2} with the previous estimates. 
		
		\begin{ethm}[Convergence of the entropic interpolation under $CD(0,n)$]\label{thm:interpolationtoflow}
			\label{gthm-3}
			Let $\mu,\nu$ be two absolutely continuous and compactly supported measures with smooth density w.r.t. $m$. For all $T>0$, $\left( \mu_t^T\right)_{t \in [0,T]}$ denotes the entropic interpolation from $\mu$ to $\nu$. Then for every $T>1$ and $t \in (0,T)$ we have
			$$
			W_2 \left(\mu_t^T,P_t^{*} \mu \right) \leq 2 \sqrt{2(\mathcal{F}(\nu)-\mathcal{F}(\mu))+\mathcal{C}_1(\mu,\nu)+2n \log(T)}\left(\sqrt{T}-\sqrt{T-t} \right).
			$$
			Furthermore for any $a>1$, there exists a constant $C>0$, such that for all $T \geq a$ and $t \in [0,a]$, 
			$$
			W_2 \left(\mu_t^T,P_t^{*} \mu \right) \leq C \sqrt{\frac{n \log(T)}{T}}.
			$$
		\end{ethm}

		\begin{erem}
			\begin{itemize}
				\item The findings of Theorem \ref{thm:interpolationtoflow} may not be optimal. More precisely, they are not optimal for $\R^n$ equipped with the usual Laplacian operator as we will see in the next section. However, we do not know whether it is possible to improve on Theorem \ref{thm:interpolationtoflow} assuming the $CD(0,n)$ condition only. The natural conjecture is that under the hypothesis of this section the convergence rate is $T^{-1}$, namely
				$$
				W_2 \left(\mu_t^T,P_t^{*} \mu \right) \leq \frac{C}{T},\,\,T>0.
				$$
				\item The $CD(0,n)$ condition is not strong enough to imply that $m$ is a probability measure: if we were to add this assumption then, combining the results of \cite{conforti-tamanini2019} and the methods of this paper, we could obtain a better convergence rate of $T^{-1/2}$.
			\end{itemize}
		\end{erem}
		
		\subsubsection{A turnpike estimate}
		
		Under the $CD(0,n)$ condition the following estimates for the Fisher information along the heat flow is well-known
		$$
		\mathcal{I}_W(P_t^*(\mu))=|\grad_{P_t^{*}\mu} \mathcal{F}|_{P_t^*(\mu)}^2 \leq \frac{n}{2t}, \ \mu \in \mathcal{P}_2(N), \, t>0.
		$$
		We can show an analogous estimate along the entropic interpolations. The proof is the exact analogous of the Theorem~\ref{gio-thm3} by using the estimate~\eqref{eq-82}.
		
		\begin{ethm}\label{thm:tpikeCD(0,N)} Let $\mu,\nu \in \mathcal{P}_2(N)$ be two compactly supported absolutely continuous measures with smooth positive densities against $m$. Then for every $T>0$ and $t \in (0,T)$ we have
			$$
			|\mathrm{grad}_{\mu_t^T} \mathcal{F}|_{\mu_t^T}^2=\mathcal{I}_W(\mu_t^T) \leq \frac{n}{{{2}}t}+\frac{n}{{{2}}(T-t)},
			$$
			that is for every $T > 0 $ and $\theta \in (0,1)$
			$$
			\mathcal{I}_W(\mu_{\theta T}^T) \leq \frac{n}{{{2}}T \theta (1- \theta)}.
			$$
		\end{ethm}
		
		\subsection{A refined study of the Euclidean heat semigroup in $\R^n$}\label{sec:heat}
		
		In this subsection $\left(P_t \right)_{t \geq 0}$ is the usual heat semigroup in $\R^n$,  in that case,  $m$ is the Lebesgue measure in $\R^n$, and  the density  kernel is given by  
		$$
		\forall x,y \in \R^n, \ \forall t>0, \ p_t(x,y)=\frac{1}{\left(4 \pi t\right)^{n/2}}e^{-\frac{|x-y|^2}{4t}}.
		$$
		Recall that  $(P_t)_{t \geq 0}$ verifies the $CD(0,n)$ curvature-dimension condition. In this setting we can improve some of the results of the former section relying on a different method that exploits $\Gamma$-convergence. The first step is to establish a $\Gamma$-convergence result analogous to the one recently proven in~\cite{conforti-tamanini2019} under the hypothesis that $m$ is a probability measure. This hypothesis is clearly violated here. For the definition and basic properties of $\Gamma$-convergence we refer to~\cite{braides2002gamma}. For $T>0$ we denote by $R_{0T}^T$ the positive measure,
		$$
		dR^T_{0T}(x,y) =  p_T(x,y) dm(x)dm(y).
		$$
		A crucial observation here is that for all $T>0$ we have $supp(f^T)=supp(\mu)$ and $supp(g^T)=supp(\nu)$. This follows from  equation~(\ref{eq-43}) at time $t=0$ and $t=T$ and the basic properties of the heat semigroup. Let us now prove the announced $\Gamma$-convergence result.
		
		
		
		
		\begin{ethm}[$\Gamma$-convergence of the Schrödinger problem] \label{thm-gammaconv}
			
			Let $\mu,\nu \in \mathcal P_2(\R^n)$ be compactly supported and absolutely continuous probability measures and $(T_{k})_{k\geq 1}$ a diverging sequence. Then the sequence of functionals 
			\begin{equation}\label{eq-176}
			H\left( \ \cdot \ | R_{0T_k}^{T_k} \right) - \frac{n}{2}\log \left( 4 \pi T_k \right) 
			\end{equation}
			defined on $\Pi(\mu,\nu)$, which we equip with the weak convergence topology, $\Gamma$-converge to the functional $ H \left( \ \cdot \ |m \otimes m\right)$. In particular,
			
			\begin{equation}\label{eq: precise cost asymptotics}
			\mathcal{C}_{T_k}(\mu,\nu)-{2n}\log(4 \pi T_k) \underset{k \rightarrow + \infty}{\rightarrow} 2\mathcal{ F}(\mu)+ 2\mathcal{ F}(\nu). 
			\end{equation}
			
		\end{ethm}
		
		After noticing that $H\left( \ \cdot \ | R_{0T_k}^{T_k} \right) - \frac{n}{2}\log \left( 4 \pi T_k \right) $ is a decreasing sequence of functionals, we could invoke \cite[Prop 5.7]{dal2012introduction} to obtain that the $\Gamma$-limit of the sequence is the loweremicontinuous envelope of the pointwise limit $H \left( \ \cdot \ |m \otimes m\right)$. Since relative entropy is lowersemicontinuous in the weak topology, this argument proves Theorem \ref{thm-gammaconv}. A direct proof can also be obtained rather easily working directly on the definition of $\Gamma$-convergence.  Therefore, in the interest of being self contained and the for the reader's convenience, we decided to include it in this manuscript.

		\begin{eproof}
			Let $(T_{k})_{k \geq 1}$ be a diverging sequence. We begin by proving the liminf inequality: consider $\gamma_k\rightarrow\gamma$ weakly and recall that 
			$$
			\frac{d \, R^{T_k}_{0T_k} }{d \,  m \otimes m}= \frac{1}{(4\pi T_k)^{n/2}} \exp\left(-\frac{|y-x|^2}{4T_k}\right),
			$$
			which gives
			$$
			H(\gamma_k | R_{0T_k}^{T_k})= H(\gamma_k | m\otimes m) + \frac{1}{4T_k}\int |x-y|^2 d\gamma_k(x,y) + \frac{n}{2}\log( 4\pi T_k)
			$$
			The desired inequality follows by letting $k\rightarrow\infty$, the lowersemicontinuity of $H(\cdot |m\otimes m)$ and the fact that for all $k$, the marginals of $\gamma_k$ are $\mu$ and $\nu$, which admits a second moment. For the limsup inequality, it suffices to choose $\gamma_k\equiv\gamma$ as recovery sequence and argue as we just did. Let us now move to the proof of \eqref{eq: precise cost asymptotics}. We first observe that the optimal coupling in 
			$$
			\underset{\gamma \in \Pi(\mu,\nu)}{\inf}H\PAR{\gamma | m \otimes m}
			$$
			is $\mu \otimes \nu$. Indeed, $\mu \otimes \nu = (\frac{d \mu}{dm}\times \frac{d \nu}{dm})(m \otimes m)$ is a transport plan between $\mu$ and $\nu$ which is also a $\PAR{f,g}$-transform of $m \otimes m$ and such transport plans are optimal in the Schrödinger problem, see \cite[Proposition 4.1.5]{tamanini2017}. Moreover it is easily checked that 
			$$
			H\PAR{\mu\otimes\nu| m \otimes m} = \mathcal{F}(\mu)+\mathcal{F}(\nu).
			$$
			Since $\Pi(\mu,\nu)$ is weakly compact, using the basic properties of $\Gamma$-convergence we have the convergence of optimal values in \eqref{eq-176}, whence \eqref{eq: precise cost asymptotics}. 
		\end{eproof}

		\begin{erem}
			In the setting of this section the expansion of $\mathcal{C}_T(\mu,\nu)-{2n}\log(4 \pi T)$ can be improved to
			$$
			T(\mathcal{C}_T(\mu,\nu)-{2n}\log(4 \pi T)-2(\mathcal{ F}(\mu)+ \mathcal{ F}(\nu))) \underset{T \rightarrow + \infty}{\rightarrow} \frac14\int |x-y|^2 d\mu \otimes\nu(x , y).
			$$
			Similar results have been obtained in \cite{feydy2019interpolating}, where the convergence of the so called Sinkhorn divergences towards MMD divergences is established.
		\end{erem}

		Let us prove the announced convergence at speed $1/T$.

		\begin{ethm}[Convergence of entropic interpolations in $\R^n$]\label{thm: precise conv to grad flow}
			Let $\mu,\nu \in \mathcal{P}_2(\R^n)$ be two compactly supported absolutely continous measures with smooth positive densities w.r.t.   $m$. If $(\mu_t^T)_{t \in [0,T]}$ is the entropic interpolation from $\mu$ to $\nu$ then for every $T>0$ and $t \in (0,T)$,
			
			\begin{equation}\label{eq-303}
			\left|(T-t)^2 |\grad_{\mu_t^T} \mathcal{F} + \dot \mu_t^T |_{\mu_t^T}^2 - \int | x- \int y d\nu(y)|^2dP_t^*\mu(x) \right| \underset{ T \rightarrow + \infty}{\rightarrow } 0.
			\end{equation}
			
			Moreover, for every $a>0$, there exists a constant $C>0$ such that for every $T>a\geq t\geq 0$,
			
			\begin{equation}\label{eq-304}
			W_2(\mu_t^T,P_t^*\mu) \leq \frac{C}{T-t}.
			\end{equation}
			
		\end{ethm}
		
		\begin{eproof}
			Let $\left(f^T,g^T \right)$ be two functions in $L^{\infty}(m)$ such that for every $t \in [0,T]$: $\mu_t^T=P_tf^T P_{T-t}g^Tdm$ and $\|g^T\|_{L^1(m)}$=1.  Observe that for all $0\leq t \leq T$
			$$
			\left|\grad_{\mu_t^T}\mathcal{ F} + \dot \mu_t^T \right|_{\mu_t^T}^2={4}\int \Gamma \PAR{\log P_{T-t}g^T}d\mu_t^T.
			$$
			Moreover $\Gamma(\log P_{T-t}g^T)=\frac{|\nabla P_{T-t}g^T|^2}{\left|P_{T-t}g^T \right|^2}=\frac{| P_{T-t}\nabla g^T|^2}{\left|P_{T-t}g^T \right|^2}$  since $\nabla P_{T-t}g^T=P_{T-t}(\nabla g^T)$ for the Euclidean heat semigroup. Hence we get, 
			\begin{equation}\label{eq: gamma g}
			\begin{split}
			2(T-t)\frac{P_{T-t}\nabla g^T(x)}{P_{T-t}g^T(x)}&= 2(T-t)\frac{\int \nabla g^T(y)p_{T-t}(x,y)dm(y)}{\int g^T(y)p_{T-t}(x,y)dm(y)} \\
			&=\frac{\int g^T(y)(x-y)p_{T-t}(x,y)dm(y)}{\int g^T(y) p_{T-t}(x,y)dm(y) }\\
			&= x - \frac{\int g^T(y)y p_{T-t}(x,y)dm(y)}{\int g^T(y) p_{T-t}(x,y)dm(y)}.\\
			\end{split}
			\end{equation}
			Next, we observe that a slight modification of Lemma 3.6\footnote{The Lemma does not apply directly since $m$ is not a probability measure. Therefore hypothesis (H2) therein is violated. However, it is not difficult to see, that in the particular case of the heat semigroup on $\R^n$, we can remove this hypothesis. We omit the details here.} in \cite{conforti-tamanini2019} yields that $g^T \rightarrow \frac{d\nu}{dm}$ in $L^2(m)$ as $T\rightarrow\infty$. Using this convergence and that $supp(g^T)=supp(\nu)$ for all $T$, we obtain   for any compact set $K\subset \R^n$, 
			$$
			\int g^T(y)e^{-\frac{|x-y|^2}{4(T-t)}}dm(y) \underset{T \rightarrow \infty}{\rightarrow} 1, \quad \text{uniformly for $x\in K$,}
			$$
			
			and
			$$
			\int g^T(y)ye^{-\frac{|x-y|^2}{4(T-t)}}dm(y) \underset{T \rightarrow \infty}{\rightarrow} \int y d\nu(y), \quad \text{uniformly for $x\in K$.}
			$$
			Therefore, if we define $\theta_T(x)=\frac{\int g^T(y)y p_{T-t}(x,y)dm(y)}{\int g^T(y) p_{T-t}(x,y)dm(y)}$, we have
			
			\begin{equation}\label{conv compact sets}
			\theta_T(x) \underset{T \rightarrow \infty}{\rightarrow} \int y d\nu(y) \quad \text{uniformly on compact sets.}
			\end{equation}
			Moreover, using the fact that $supp(g^T)=supp(\nu)$, there exists a constant $C$ such that 
			\begin{equation}\label{theta bound}
			\forall x\in \R^n,T\geq1,\quad |\theta_T(x)|\leq C.
			\end{equation}
			Therefore we have 
			\begin{equation*}
			\begin{split}
			{}&\left|4(T-t)^2 \int \Gamma \PAR{\log (P_{T-t}g^T)}(x)d\mu_t^T(x)- \int |x- \int y d \nu(y)|^2dP_t^*\mu(x) \right| \\ 
			&=\left| \int |x-\theta_T(x)|^2 d\mu_t^T(x)- \int |x- \int y d \nu(y)|^2dP_t^*\mu(x) \right|\\
			&\leq\left| \int |x|^2d\mu_t^T(x)-\int |x|^2 d P_t^*\mu(x)   \right|+\left| \int |\theta_T(x)|^2 d\mu^T_t(x) - \left|\int y d\nu(y) \right|^2  \right|\\
			& +\left|\int 2\langle x,\theta_T(x)\rangle d\mu^T_t(x) - \int 2\langle x,\int y d \nu(y)\rangle d P^*_t\mu(x)  \right|.
			\end{split}
			\end{equation*}
			Since $\mu_t^T \overset{W_2}{\underset{T \rightarrow \infty}{\rightarrow}} P_t^* \mu$ by Theorem \ref{thm:interpolationtoflow}, the first term in the above display vanishes as $T\rightarrow+\infty$. Using \eqref{conv compact sets},\eqref{theta bound} and the the fact the second moment of $\mu^T_t$ is uniformly bounded in $T$, we also obtain that the second term vanishes. The third term is bounded above by 
			\begin{equation*}
			\left|\int 2\langle x,\theta_T(x)-\int y d \nu(y)\rangle d\mu^T_t(x) \right|+ \left|\int 2\langle x,\int y d \nu(y) \rangle d\mu^T_t(x) - \int 2\langle x,\int y d \nu(y)\rangle d P^*_t\mu(x)  \right|.
			\end{equation*}
			Using again $\mu_t^T \overset{W_2}{\underset{T \rightarrow \infty}{\rightarrow}} P_t^* \mu$ we get,
			\begin{equation*}
			\left|\int 2\langle x,\int y d \nu(y) \rangle d\mu^T_t(x) - \int 2\langle x,\int y d \nu(y)\rangle d P^*_t\mu(x)  \right|\underset{T \rightarrow \infty}{\rightarrow}0.
			\end{equation*}
			Moreover, for all  $M>0$ fixed we have from \eqref{conv compact sets} that
			\begin{equation*}\label{eq-301}
			\left|\int_{\{|x|\leq M\}}  \langle x,\theta_T(x)-\int y d\nu(y)\rangle d\mu^T_t(x) \right|\underset{T \rightarrow \infty}{\rightarrow}0.
			\end{equation*}
			Moreover, by Cauchy Schwartz, \eqref{theta bound} and Markov's inequality,
			\begin{equation*}
			\begin{split}
			\left|\int_{\{|x|\geq M\}}  \langle x,\theta_T(x)-\int y d\nu(y)\rangle d\mu^T_t(x) \right|
			&\leq \frac{2 C}{M} \int |x|^2 d\mu^T_t(x).
			\end{split}
			\end{equation*}
			Finally, observe that $\int |x|^2 d\mu^T_t(x)$ is uniformly bounded in $T$ by a constant $D$, we have obtained 
			$$
			\forall M>0,\quad \limsup_{T\rightarrow +\infty} \left| \int \langle x,\theta_T(x)-\int y d\nu(y)\rangle d\mu^T_t(x) \right|\leq \frac{2CD}{M},
			$$
			from which the claim \eqref{eq-303} follows. The remaining claim \eqref{eq-304} is obtained by repeating the proof of Theorem \ref{thm:interpolationtoflow} replacing \eqref{eq-305} with the stronger bound \eqref{eq-303}.
		\end{eproof}
		
		\black
		
		\appendix
		
		\section{Details about the examples} \label{AppendiceA}

		\subsection{A $(0,n)$-convex function}\label{AppendiceA2}
		
		To understand what happened in the $(0,n)$-convex case, let's begin by an example on the real line. The prototypal $(0,1)$-convex function is $F(x)=-\log x$, $x>0$.  This is a $(0,1)$-convex function since
		$
		F''=(F')^2.
		$
		Let $x,T>0$, for simplicity we just treat the case where $x=y$. The gradient flow from $x>0$, denoted by $(S_t(x))_{t\geq0}$ is the solution of the ODE 
		$\dot X_t={1}/{X_t}$ starting from $x$, 
		hence for all $t>0$, $S_t(x)=\sqrt{2t+x^2}$. 
		
		The Newton system associated  is $$
		\left\{
		\begin{array}{l}
		\ddot X_t=-\frac{1}{X_t^3}, \\
		X_0=X_T=x.
		\end{array}
		\right.
		$$
		Now $(X_t^T)_{t \in [0,T]}$ denote the entropic interpolation between $x$ and $x$. The conserved quantity is given by 
		$
		E_T(x,x)=\dot {X_t^T}^2-\frac{1}{{X_t^T}^2}.
		$
		Thus $\left|\dot X_t^T\right|=\sqrt{E_T(x,x)+F'(X_t^T)^2}$ and we can deduce that
		$$
		\left\{\begin{array}{cc}
		\dot X_t^T=\sqrt{E_T(x,x)+F'(X_t^T)^2},& t \in [0,T/2]; \\
		\dot X_t^T = - \sqrt{E_T(x,x)+F'(X_t^T)^2},& t \in (T/2,T]. 
		\end{array} \right.
		$$
		In this example we have enough information to compute explicitly the conserved quantity.
		\begin{eprop} For $T>0$ and $x \in \R$,
			$
			E_T(x,x)= 2\frac{-x^2- \sqrt{x^4+T^2}}{T^2}.
			$
		\end{eprop}

		\begin{eproof}
			By the continuity in $T/2$ of $\left(X_t^T \right)_{t \in [0,T ]}$ we can deduce that $E_T(x,x)=-F'(X_{T/2}^T)^2$. Notice that for all $t \in [0,T/2]$
			$$
			\frac{\dot X_t^T}{F'(X_t^T)}= \sqrt{1+ \frac{E_T(x,x)}{F'(X_t^T)^2}}
			$$
			and 
			$$
			\frac{d}{dt} \frac{\dot X_t^T}{F'(X_t^T)}=-\frac{F''(X_t^T)}{F'(X_t^T)^2}\mathcal{E}_T(x,x)=-\mathcal{E}_T(x,x).
			$$
			By integration of this inequality we see that for every $t \in [0,T/2)$
			\begin{equation} 
			\label{geq-2}
			\sqrt{1+\frac{E_T(x,x)}{F'(X_t^T)^2}}-\sqrt{1+\frac{E_T(x,x)}{F'(x)^2}}=\frac{TE_T(x,x)}{2}. 
			\end{equation}
			When  $t=T/2$ we get 
			$
			\frac{T^2}{4}E_T(x,x)^2- \frac{E_T(x,x)}{F'(x)^2}-1=0
			$
			and since  $E_T(x,x) \leq 0$ we deduce that
			$$
			E_T(x,x)=\frac{- \frac{1}{F'(x)^2}-\sqrt{\frac{1}{F'(x)^4}+T^2}}{T^2/2}=\frac{-x^2-\sqrt{x^4+T^2}}{T^2/2}.
			$$
		\end{eproof}
		Hence $E_T(x,x)$ is of order $1/T$ in this  case. From~\eqref{geq-2}, we can deduce an explicit formula for $X_t^T$.
		
		\begin{eprop}
			For $x \in \R$ and $T>0$, the entropic interpolation from $x$ to $x$ is given by
			$$
			X_t^T= \sqrt{x^2+t^2 E_T(x,x)+2t\sqrt{1+E_T(x,x)x^2}}, \, 0 \leq t \leq T.
			$$
			Furthermore $X_t^T \rightarrow S_t(x)$ when $T \rightarrow \infty$, more precisely  
			$$
			X_t^T-S_t(x) \underset{T \rightarrow \infty}{\sim} E_T(x,x)\frac{t+t^2x^2}{2\sqrt{x^2+2t}},
			$$
			hence there exists a constant $C>0$ such that 
			$
			|X_t^T-S_t(x)|  \underset{T \rightarrow \infty}{\sim} \frac{C}{T}. 
			$
		\end{eprop}
		In this particular case we can compute the cost in an explicit way. 
		\begin{eprop}
			For every $x \in \R$,
			$
			C_T(x,x)  \underset{T \rightarrow \infty}{\sim} 2 \log(T).
			$
		\end{eprop}
		
		\begin{eproof}
			By the very defnition of the cost,
			\begin{equation*}
			\begin{split}
			C_T(x,x)&=\int_0^T \left(\dot {X_t^T}^2+\frac{1}{{X_t^T}^2} \right)dt=2\int_0^{T/2} \left(\dot {X_t^T}^2+\frac{1}{{X_t^T}^2} \right)dt \\
			&=4 \int_0^{T/2} \dot {X_t^T}^2dt+2 \int_0^{T/2} \left(\frac{1}{X_t^T}- \dot {X_t^T}^2 \right)dt \\
			&=4 \int_0^{T/2} \dot X_t^T \frac{\sqrt{1+{X_t^T}^2E_T(x,x)}}{X_t^T}-T E_T(x,x) \\
			&=\int_{\sqrt{-E_T(x,x)}x}^1 \frac{\sqrt{1-v^2}}{v}dv -T E_T(x,x). \\
			\end{split}
			\end{equation*}
			Hence,
			$
			C_T(x,x) \underset{T \rightarrow \infty}{\sim} 2 \log (T). 
			$
		\end{eproof}

		\subsection{The example of two gaussians on $\R$}\label{AppendiceA1}
		
		This example take place on $\R$. This is a flat space of dimension one, that mean it verify the $CD(0,1)$ condition. Recall that for $m \in \R$ and $\sigma>0$ the normal law of expected value $m$ and variance $\sigma^2$ is the probability measure on $\R$ with density against the Lebesgue measure, 
		$$
		\mathcal{N}(m,\sigma^2)(x)=\frac{1}{\sqrt{2\pi \sigma^2}}\exp \left(- \frac{(x-m)^2}{2 \sigma^2} \right). 
		$$
		In this case we know an expression for the heat semigroup, for $f \in L^{\infty}(\R)$, we have 
		$$
		\forall t>0, \ P_tf=\mathcal{N}(0,2t)*f.
		$$
		Furthermore, for all $m \in \R$ and $\sigma,t>0$ we know that from elementary probability theory
		$$
		\mathcal{N}(0,2t)*\mathcal{N}(m,\sigma^2)=\mathcal{N}(m,\sigma^2+2t). 
		$$
		Thanks to all of these considerations we can make explicit calculus in this case. For simplicity here we are gonna consider the case of two centered gaussian measure, that is $\sigma^2=1$. Let $x_0,x_1 \in \R$, $T>0$, $\mu=\mathcal{N}(x_0,1)$ and $\nu = \mathcal{N}(x_1,1)$. We can solve explicitely the Schrödinger system
		$$
		\left\{ \begin{array}{cc}
		\mu=fP_Tg , \\
		\nu=gP_Tf,
		\end{array} \right.
		$$
		by searching solutions of the form $x \mapsto a \exp\left(- \frac{(x-b)^2}{2c^2} \right)$ with $a,b,c \in \R$. We can make explicit computations to find two solutions given by for all $x \in \R$
		$$
		\left\{ \begin{array}{cc}
		f(x)=\dfrac{1}{\sqrt{2 \pi \mathcal{D}_T^2}} \exp \left(- \frac{\left(x - \frac{\mathcal{D}_T^2(\mathcal{D}_T^2+2T)^2}{(\mathcal{D}_T^2+2T)^2-\mathcal{D}_T^2} \left(x_0-\frac{\mathcal{D}_T^2}{\mathcal{D}_T^2+2T}x_1 \right) \right)^2}{2\mathcal{D}_T^2}\right), \\
		g(x)= \sqrt{\mathcal{D}_T^2+2T} \exp \left(- \frac{\left(x - \frac{\mathcal{D}_T^2(\mathcal{D}_T^2+2T)^2}{(\mathcal{D}_T^2+2T)^2-\mathcal{D}_T^2} \left(x_1-\frac{\mathcal{D}_T^2}{\mathcal{D}_T^2+2T}x_0 \right)\right)^2}{2\mathcal{D}_T^2}\right), 
		\end{array} \right.
		$$
		where the parameter $\mathcal{D}_T$ is given by
		$
		\mathcal{D}_T^2=\sqrt{(T-1)^2+2T}-(T-1).
		$
		Observe that $f$ is the density of the normal law $\mathcal{N}\left(\frac{\mathcal{D}_T^2(\mathcal{D}_T^2+2T)^2}{(\mathcal{D}_T^2+2T)^2-\mathcal{D}_T^2} \left(x_0-\frac{\mathcal{D}_T^2}{\mathcal{D}_T^2+2T}x_1 \right),\mathcal{D}_T^2 \right)$. This is an arbitrary choice, because there is only unicity up to the trivial transform $(f,g) \mapsto (cf,g/c)$ for some $c \in \R$. From those expressions we can easily deduce a formula for the entropic interpolation $(\mu_t^T)_{t \in [0,T]}$ between $\mu$ and $\nu$, actually it's a normal law $\mathcal{N}(x_t^T,\sigma_t^T)$ where the parameter are given by
		\begin{equation*}
		\left\{\begin{array}{cc}
		x_t^T=\frac{T-t}{T}x_0+\frac{t}{T}x_1,  \\
		\sigma_t^T=1+2\frac{t(T-t)}{\mathcal{D}_T^2+T}.  
		\end{array} \right.
		\end{equation*}
		We want to quantify the convergence of $\mu^T$ toward the gradient flow $\left(P_t^* (\mu) \right)_{t \in [0,T]}$. The gradient flow is given by $P_t^*(\mu)=P_t\left(\frac{d\mu}{dm} \right)dm=\mathcal{N}(x_0,\mathcal{D}_T^2+2t)$. Actually the Wasserstein distance between two gaussian measures can be explicitely computed. Indeed let $\mu= \mathcal{N}(m_0,\sigma_0^2)$ and $\nu=\mathcal{N}(m_1,\sigma_1^2)$, the map $T:x \mapsto \frac{\sigma_1}{\sigma_0}(x-m_0)+m_1$ verify $T\#\mu=\nu$, hence by Brenier theorem $W_2^2(\mu,\nu)=\int |x-T(x)|^2 d \mu(x)$. From this expression and some easy computations we find
		$$
		W_2^2(\mu,\nu)=|\sigma_0-\sigma_1|^2+|m_0-m_1|^2.
		$$
		For the detail and the extension to Gaussian vectors we refer to \cite[Remark 2.31]{COTFNT}. Hence we can compute explicitly the Wasserstein distance between the entropic interpolation and the gradient flow.
		
		\begin{eprop} In the notations of this subsection
			$$
			W_2^2\left(\mu_t^T,P_t^* \mu  \right)=\frac{t^2}{T^2}(x_0-x_1)^2+\left|\sqrt{\sigma_t^T}-\sqrt{2t+1} \right|^2,
			$$
			and there exists a constant $C>0$ such that
			$
			W_2^2(\mu_t^T,P_t^* \mu) \underset{T \rightarrow \infty}{\sim} \frac{C}{T^2}.
			$
		\end{eprop}
		The velocity of $\left( \mu_t^T \right)_{t \in [0,T]}$ is given by 
		$$
		\forall t\in[0,T], \ \forall x \in \R, \ \dot \mu_t^T(x)=  \frac{\dot \sigma_t^T}{2\sigma_t^T}\left( x-x_{t/T}\right)+ \frac{1}{T} \left(x_1 - x_0 \right).
		$$
		Now we have all the element we need to compute the conserved quantity, and the following proposition follow from basic integration. 
		\begin{eprop} In the notations of this subsection we have the following equality for every $T>0$,
			$$
			\mathcal{E}_T(\mu,\nu):={\left| \dot \mu_t^T \right|_{{\mu_t}^T}^2}- {\left|\n \log \left(\mu_t^T \right) \right|_{{\mu_t}^T}^2}=\frac{\dot {\sigma_t^T}^2}{4 \sigma_t^T}+\frac{1}{T^2}(x_1-x_0)^2 - \frac{1}{\sigma_t^T}, \, 0 \leq t \leq T.
			$$
			In particular we can take $t=T/2$ to find,
			$
			\mathcal{E}_T(\mu,\nu) \underset{T \rightarrow \infty}{\sim}\frac{(x_1-x_0)^2}{T^2}-\frac{2}{T+2}
			$
			and finally we get
			$$
			\mathcal{C}_T(\mu,\nu) \underset{T \rightarrow \infty}{\sim} 2 \log (T).
			$$
		\end{eprop}
		
		\medskip
		
		\paragraph{Acknowledgements} This research was supported by the French ANR-17-CE40-0030 EFI project. We would also like to warmly
		thank the referee who gave us useful comments and pointed out references.
		
		\footnotesize
		

\begin{thebibliography}{BVCGL20}
			
			\bibitem[AGS08]{ambrosio-gigli2008}
			L.~{Ambrosio}, N.~{Gigli}, and G.~{Savar\'e}.
			\newblock {\em {Gradient flows in metric spaces and in the space of probability
					measures. 2nd ed.}}
			\newblock Basel: Birkh\"auser, 2nd ed. edition, 2008.
			
			\bibitem[BB00]{benamou-brenier2000}
			J.-D. {Benamou} and Y.~{Brenier}.
			\newblock {A computational fluid mechanics solution to the Monge-Kantorovich
				mass transfer problem.}
			\newblock {\em {Numer. Math.}}, 84(3):375--393, 2000.
			
			\bibitem[BE85]{bakry-emery1985}
			D.~{Bakry} and M.~{Emery}.
			\newblock {Diffusions hypercontractives.}
			\newblock {\em {S\'emin. de probabilit\'es XIX, Univ. Strasbourg 1983/84,
					Proc., Lect. Notes Math. 1123, 177-206}}, 1985.
			
			\bibitem[BGL14]{bgl-book}
			D.~{Bakry}, I.~{Gentil}, and M.~{Ledoux}.
			\newblock {\em {Analysis and geometry of Markov diffusion operators.}}
			\newblock Cham: Springer, 2014.
			
			\bibitem[BL06]{bakry-ledoux}
			D.~{Bakry} and M.~{Ledoux}.
			\newblock {A logarithmic Sobolev form of the Li-Yau parabolic inequality.}
			\newblock {\em {Rev. Mat. Iberoam.}}, 22(2):683--702, 2006.
			
			\bibitem[{Bra}02]{braides2002gamma}
			A.~{Braides}.
			\newblock {\em {\(\Gamma\)-convergence for beginners.}}, volume~22.
			\newblock Oxford University Press, 2002.
			
			\bibitem[BVCGL20]{backhoff2019mean}
			J.~Backhoff-Veraguas, G.~Conforti, I.~Gentil, and C.~L{\'e}onard.
			\newblock The mean field schr{\"o}dinger problem: ergodic behavior, entropy
			estimates and functional inequalities.
			\newblock {\em Probab. Theory Related Fields}, 178(1-2):475--530, 2020.
			
			\bibitem[CCG21]{ccg2}
			G.~Clerc, G.~Conforti, and I.~Gentil.
			\newblock On the variational interpretation of local logarithmic sobolev
			inequalities, Preprint 2021.
			
			\bibitem[CGP16]{chengeorgiou2016}
			Y.~Chen, T.~T. Georgiou, and M.~Pavon.
			\newblock On the relation between optimal transport and {S}chr\"odinger
			bridges: a stochastic control viewpoint.
			\newblock {\em J. Optim. Theory Appl.}, 169(2):671--691, 2016.
			
			\bibitem[Cle21]{clerc}
			G~Clerc.
			\newblock {Regularity of the Schr{\"o}dinger cost.}
			\newblock Preprint, 2021.
			
			\bibitem[Con19]{conforti2017}
			G.~Conforti.
			\newblock {A second order equation for Schr\"odinger bridges with applications
				to the hot gas experiment and entropic transportation cost.}
			\newblock {\em {Probab. Theory Relat. Fields}}, 174(1-2):1--47, 2019.
			
			\bibitem[{Cos}85]{costa}
			M.~H.~M. {Costa}.
			\newblock {A new entropy power inequality.}
			\newblock {\em {IEEE Trans. Inf. Theory}}, 31:751--760, 1985.
			
			\bibitem[CT21]{conforti-tamanini2019}
			G.~{Conforti} and L.~{Tamanini}.
			\newblock {A formula for the time derivative of the entropic cost and
				applications}.
			\newblock {\em {J. Funct. Anal.}}, 280(11):48, 2021.
			
			\bibitem[DM12]{dal2012introduction}
			G.~Dal~Maso.
			\newblock {\em An introduction to $\Gamma$-convergence}, volume~8.
			\newblock Springer Science \& Business Media, 2012.
			
			\bibitem[EKS15]{erbar-kuwada2015}
			M.~{Erbar}, K.~{Kuwada}, and K.-T. {Sturm}.
			\newblock {On the equivalence of the entropic curvature-dimension condition and
				Bochner's inequality on metric measure spaces.}
			\newblock {\em {Invent. Math.}}, 201(3):993--1071, 2015.
			
			\bibitem[{Erb}10]{erbar}
			M.~{Erbar}.
			\newblock {The heat equation on manifolds as a gradient flow in the Wasserstein
				space}.
			\newblock {\em {Ann. Inst. Henri Poincar\'e, Probab. Stat.}}, 46(1):1--23,
			2010.
			
			\bibitem[{Fol}88]{follmer1988}
			H.~{Follmer}.
			\newblock {Random fields and diffusion processes.}
			\newblock {Calcul des probabilit\'es, \'Ec. d'\'Et\'e, Saint-Flour. 1985-87,
				Lect. Notes Math. 1362, 101-203}, 1988.
			
			\bibitem[FSV{\etalchar{+}}19]{feydy2019interpolating}
			J.~Feydy, T.~S{\'e}journ{\'e}, F.-X. Vialard, S.-I. Amari, A.~Trouv{\'e}, and
			G.~Peyr{\'e}.
			\newblock Interpolating between optimal transport and mmd using sinkhorn
			divergences.
			\newblock In {\em The 22nd International Conference on Artificial Intelligence
				and Statistics}, pages 2681--2690, 2019.
			
			\bibitem[{Gig}12]{gigli2012}
			N.~{Gigli}.
			\newblock {Second order analysis on $(\mathcal P_{2}(M),W_{2})$.}
			\newblock {\em {Mem. Am. Math. Soc.}}, 1018:154, 2012.
			
			\bibitem[GLR17]{gentil-leonard2017}
			I.~Gentil, C.~L\'eonard, and L.~Ripani.
			\newblock About the analogy between optimal transport and minimal entropy.
			\newblock {\em Ann. Fac. Sci. Toulouse Math. (6)}, 26(3):569--601, 2017.
			
			\bibitem[GLR20]{gentil-leonard2020}
			I.~{Gentil}, C.~{L\'eonard}, and L.~{Ripani}.
			\newblock {Dynamical aspects of the generalized Schr\"odinger problem via Otto
				calculus -- a heuristic point of view}.
			\newblock {\em {Rev. Mat. Iberoam.}}, 36(4):1071--1112, 2020.
			
			\bibitem[GT20]{gigli-tamanini2020}
			N.~{Gigli} and L.~{Tamanini}.
			\newblock { Benamou-Brenier and Kantorovich duality formulas for the entropic
				cost on $RCD^*(K,N)$ spaces }.
			\newblock {\em {Probab. Theory Relat. Fields}}, 176(1-2):1--34, 2020.
			
			\bibitem[GT21]{gigli2018second}
			N.~Gigli and L.~Tamanini.
			\newblock Second order differentiation formula on {RCD}*{$(K,N)$} spaces.
			\newblock {\em J. Eur. Math. Soc. (JEMS)}, 23(5):1727--1795, 2021.
			
			\bibitem[JKO98]{jko1998}
			R.~{Jordan}, D.~{Kinderlehrer}, and F.~{Otto}.
			\newblock {The variational formulation of the Fokker-Planck equation.}
			\newblock {\em {SIAM J. Math. Anal.}}, 29(1):1--17, 1998.
			
			\bibitem[LB91]{lafrance1991envelope}
			J.~T LaFrance and L~D. Barney.
			\newblock The envelope theorem in dynamic optimization.
			\newblock {\em Journal of Economic Dynamics and Control}, 15(2):355--385, 1991.
			
			\bibitem[Lé14]{leonard2014}
			C.~Léonard.
			\newblock {A survey of the Schr\"odinger problem and some of its connections
				with optimal transport.}
			\newblock {\em {Discrete Contin. Dyn. Syst.}}, 34(4):1533--1574, 2014.
			
			\bibitem[McK63]{mckenzie1963turnpike}
			L.~McKenzie.
			\newblock Turnpike theorems for a generalized {L}eontief model.
			\newblock {\em Econometrica (pre-1986)}, 31(1, 2):165, 1963.
			
			\bibitem[{Mik}04]{mikami04}
			T.~{Mikami}.
			\newblock {Monge's problem with a quadratic cost by the zero-noise limit of
				$h$-path processes.}
			\newblock {\em {Probab. Theory Relat. Fields}}, 129(2):245--260, 2004.
			
			\bibitem[{Ott}01]{otto2001}
			F.~{Otto}.
			\newblock {The geometry of dissipative evolution equations: The porous medium
				equation.}
			\newblock {\em {Commun. Partial Differ. Equations}}, 26(1-2):101--174, 2001.
			
			\bibitem[OV00]{otto-villani2000}
			F.~{Otto} and C.~{Villani}.
			\newblock {Generalization of an inequality by Talagrand and links with the
				logarithmic Sobolev inequality.}
			\newblock {\em {J. Funct. Anal.}}, 173(2):361--400, 2000.
			
			\bibitem[PC19]{COTFNT}
			G.~Peyr{é} and M.~Cuturi.
			\newblock Computational optimal transport.
			\newblock {\em Foundations and Trends in Machine Learning}, 11(5-6):355--607,
			2019.
			
			\bibitem[{Rip}19]{ripani}
			L.~{Ripani}.
			\newblock {Convexity and regularity properties for entropic interpolations.}
			\newblock {\em {J. Funct. Anal.}}, 277(2):368--391, 2019.
			
			\bibitem[Sch31]{Sch31}
			E.~Schr\"odinger.
			\newblock {\"U}ber die {U}mkehrung der {N}aturgesetze.
			\newblock {\em Sitzungsberichte Preuss. Akad. Wiss. Berlin. Phys. Math.},
			144:144--153, 1931.
			
			\bibitem[Sch32]{Sch32}
			E.~Schr{\"o}dinger.
			\newblock Sur la th{\'e}orie relativiste de l'{\'e}lectron et
			l'interpr{\'e}tation de la m{\'e}canique quantique.
			\newblock {\em Ann. Inst. H. Poincar\'e}, 2:269--310, 1932.
			\newblock 
			
			\bibitem[{Tam}17]{tamanini2017}
			L.~{Tamanini}.
			\newblock {Analysis and geometry of RCD spaces via the Schr{\"o}dinger
				problem}.
			\newblock {\em Phd Thesis}, 2017.
			
			\bibitem[TZ15]{trelat2015turnpike}
			E.~Tr{\'e}lat and E.~Zuazua.
			\newblock The turnpike property in finite-dimensional nonlinear optimal
			control.
			\newblock {\em Journal of Differential Equations}, 258(1):81--114, 2015.
			
			\bibitem[{Vil}09]{villani2009}
			C.~{Villani}.
			\newblock {\em {Optimal transport. Old and new.}}
			\newblock Berlin: Springer, 2009.
			
			\bibitem[vS05]{renesse-sturm2005}
			M.-K. {von Renesse} and K.-T. {Sturm}.
			\newblock {Transport inequalities, gradient estimates, entropy and Ricci
				curvature.}
			\newblock {\em {Commun. Pure Appl. Math.}}, 58(7):923--940, 2005.
			
		\end{thebibliography}
		
		\newcommand{\etalchar}[1]{$^{#1}$}

	\end{document}